\documentclass[11pt,letterpaper]{article}

\usepackage[T1]{fontenc}
\usepackage{lmodern}
\usepackage{import}

\usepackage[margin=1in,marginpar=0.8in]{geometry}
\usepackage{amsmath}
\usepackage{mathtools,amssymb,amstext,mathrsfs}
\usepackage{amsthm}
\usepackage{authblk}
\usepackage{enumerate}
\usepackage[shortlabels]{enumitem}
\usepackage{tikz}
\usepackage{thmtools}
\usepackage{thm-restate}
\usepackage[normalem]{ulem}

\usepackage{float,subcaption}
\usepackage{bookmark}
\usepackage{lineno}
\usepackage{cancel}
\newtheorem{thm}{Theorem}[section]

\newtheorem{lemma}[thm]{Lemma}

\newtheorem{obs}[thm]{Observation}

\numberwithin{equation}{section}
\theoremstyle{definition} 
\newtheorem{definition}[thm]{Definition}

\renewcommand{\subset}{\subseteq}

\newcommand{\del}{\mathsf{Del}}
\newcommand{\delsave}{\mathsf{DelSave}}
\newcommand{\Omit}{\mathsf{Omit}}

\definecolor{asparagus}{rgb}{0.53, 0.66, 0.42}

\newcommand\CCC{7}
\definecolor{cerulean}{rgb}{0.0, 0.48, 0.65}
\definecolor{cornellred}{rgb}{0.7, 0.11, 0.11}
\definecolor{darklavender}{rgb}{0.45, 0.31, 0.59}
\definecolor{darkslateblue}{rgb}{0.28, 0.24, 0.55}
\definecolor{burntorange}{rgb}{0.8, 0.33, 0.0}
\date{}

\definecolor{nicelavender}{RGB}{153, 128, 250}

\usepackage[obeyclassoptions=true,subpreambles=true]{standalone}

\begin{document}


\title{Weak Degeneracy with Availability for Planar Graphs}
\author{Zdeněk Dvořák \thanks{Charles University, Prague, Czechia, email: \texttt{rakdver@iuuk.mff.cuni.cz}.  Supported by ERC-CZ project LL2328 (Beyond the Four Color Theorem) of the Ministry of Education of Czech Republic.}  \and
        Evelyne Smith-Roberge \thanks{Illinois State University, Illinois, USA, email: \texttt{esmithr@ilstu.edu}.}}
\date{\today}

\maketitle
\begin{abstract}
Every planar graph admits exponentially many proper colorings for any assignment of lists of five allowed colors to its vertices. In this paper, we show that very little beyond degeneracy is needed to prove that this statement holds. In particular, given a planar graph $G$ (but not the list assignment), our argument provides an algorithm consisting entirely of instructions of the following form:
``Color the vertex $v$ using any color $c$ from its list, and remove $c$ from the lists of neighbors of $v$.'' and ``For the edge $uv$, choose any $k$ colors from the list of $u$, then color $v$ by any color $c$ from its list that does not appear among these
$k$ colors, and remove $c$ from the lists of neighbors of $v$.''
For every assignment of lists of size five to vertices of $G$, following this algorithm in all possible ways guarantees that at least $2^{\frac{|V(G)|}{7}}$ different proper colorings of $G$ are produced (and at least one color choice is guaranteed to be available in each step).
We also show that this result has the following consequence: for every planar graph $G$, there exists a set $R\subseteq V(G)$ of size at least $0.06|V(G)|$ such that
$G$ can be properly colored from any assignment of lists of size four to vertices in $R$ and lists of size five to the rest of vertices of $G$. Moreover, these statements remain true even in stronger settings, such as for correspondence coloring.
\end{abstract}

\section{Introduction}\label{sec:intro}

In 1994, Thomassen~\cite{thomassen5LC} famously provided an elegant proof that every planar graph is $5$-list-colorable\footnote{%
Throughout the paper, all graph colorings are required to be \emph{proper}, i.e., to assign different colors to adjacent vertices, unless specified otherwise.
Let $L$ be an assignment $L$ of lists of colors to vertices of a graph $G$.  A coloring $\varphi$ of $G$
is an \emph{$L$-coloring} if $\varphi(v)\in L(v)$ holds for every vertex $v\in V(G)$, and the graph $G$ is \emph{$L$-colorable}
if an $L$-coloring of $G$ exists.  For an integer $k$, the list assignment $L$ is a \emph{$k$-list-assignment} if $|L(v)|\ge k$ holds
for every vertex $v\in V(G)$.  We say that a graph $G$ is \emph{$k$-list-colorable} if it is $L$-colorable for every $k$-list-assignment $L$ for $G$.}.
Given the simplicity of the argument, it is natural to ask whether this result can be pushed further; and indeed, we now know many strengthenings and generalizations.  For example, projective planar graphs are 5-list-colorable if and only if their clique number is at most five~\cite{bohmelc},
and more generally, the 5-list-colorability of graphs on any fixed surface is characterized by finitely many forbidden subgraphs~\cite{lukethe}.
Graphs drawn in the plane with crossings far apart are 5-list-colorable~\cite{dvovrak20175}.  When 5-list-coloring planar graphs, the colors of arbitrarily many
sufficiently distant vertices can be prescribed~\cite{precolver}.

Moreover, we know that many of these results remain true when we consider more general forms of coloring.
\begin{itemize}
\item Suppose we are not given the list assignment all at once.  Instead, an adversary reveals for each color the set of vertices whose lists
contain it, and each time forces us to decide for which of them we are going to use this color before revealing the next one.
This leads to the notion of \emph{online list colorability} (or \emph{paintability}).  Even in this stronger setting, planar graphs
are known to be online 5-list-colorable~\cite{schauz2009mr}!
\item As one of the first instances of the polynomial method, Alon and Tarsi~\cite{alon1992colorings} showed that if certain monomials have non-zero coefficients in a particular polynomial associated to the graph,
then the graph is $k$-list-colorable.  This motivated the definition of the \emph{Alon-Tarsi number} of the graph as the minimum $k$ for which this argument works.
Thus, a graph of Alon-Tarsi number $k$ is $k$-list-colorable, and actually even online $k$-list-colorable~\cite{schauz2010flexible}.
Planar graphs are known to have Alon-Tarsi number at most five~\cite{atar}, strengthening the aforementioned online list colorability result.
\item Another popular variation on list coloring (incomparable with the Alon-Tarsi number) is \emph{correspondence coloring}, where the notion of proper coloring is modified.  More precisely,
in addition to the list assignment $L$, we are also given a \emph{correspondence assignment}, which specifies for each edge $uv$ a matching between $L(u)$ and $L(v)$.
In correspondence coloring, we are then forbidden to use the matched colors at $u$ and $v$ (the usual proper coloring is the special case where each
color is matched with itself).  When Dvořák and Postle~\cite{dvovrak2018correspondence} introduced this variant of graph coloring, they noted that
Thomassen's 5-list-colorability proof works in this setting as well, and thus planar graphs are $5$-correspondence-colorable.
\end{itemize}

In all these settings, we can also ask whether the graph can indeed only have the single coloring guaranteed by the aforementioned results,
or whether it has to have many.  In general, it turns out the latter is the case, and there actually are exponentially many different colorings.
Thomassen~\cite{thomassen2007exponentially} proved this to be the case for 5-list coloring of planar graphs.
\begin{thm}[Thomassen~\cite{thomassen2007exponentially}]\label{thm:expmany}
Every $n$-vertex planar graph has at least $2^{n/9}\ge 1.08^n$ different colorings from any assignment of lists of size five to its vertices.
\end{thm}
Bosek et al.~\cite{bosek2022graph} provided a simpler proof using the polynomial method, which also improves the bound to $5^{n/4}\ge 1.49^n$.
Postle and Thomas~\cite{PosThoHyperb} extended the result to \emph{locally planar graphs}, that is graphs drawn on a surface so that
each non-contractible cycle is long (as a function of the genus of the surface).  Furthermore, Postle and Smith-Roberge~\cite{PS23}
generalized these results to the correspondence coloring setting.

Although the notions of (online) list coloring, correspondence coloring, and Alon-Tarsi number are closely related, the results we
discussed above are proved by a wide range of techniques, often specifically tailored to the particular setting.  Recently, there
has been interest in finding the simplest common principles that would be applicable in all of the settings yet still be sufficiently
powerful to imply (at least) these fundamental results.

A basic principle of this kind is \emph{degeneracy}, which can be framed as follows.
Let $G$ be a graph and let $f:V(G)\to\mathbb{Z}$ be a function.  We view this function as giving a strict lower bound on the size of the list of each vertex;
i.e., imagine we also have a list assignment $L$ such that $|L(v)|>f(v)$ holds for every vertex $v\in V(G)$.
Let $\del(G, f, v) = (G-v, f')$, where $f':V(G-v) \to \mathbb{Z}$ is given by 
$$f'(u) =\begin{cases*} 
        f(u)-1 &if $u\in N_G(v)$ \\
        f(u) &otherwise.
        \end{cases*}
$$

We say that the $\del(G,f,v)$ operation is \emph{legal} if $f(v)\ge 0$.  Intuitively, this operation reflects the state after we select an arbitrary color $c$ from the list of $v$
(which exists when the operation is legal, since then $|L(v)|>f(v)\ge 0$), use $c$ to color the vertex $v$ (removing the vertex $v$ from further consideration), and remove the color $c$
from the lists of neighbors of $v$ (here we are pessimistic and consider the worst-case situation where all these lists actually contain $c$).
If there exists an ordering $v_1$, \ldots, $v_n$ of the vertices of $G$ so that the sequence of operations $\del(G,f,v_1)=(G_1,f_1)$, $\del(G_1,f_1,v_2)=(G_2,f_2)$, \ldots, $\del(G_{n-1},f_{n-1},v_n)$
is legal, then we say that the graph $G$ is \emph{$f$-degenerate}.  If this is the case for the function that assigns the same value $d$ to each vertex, then we say that the graph is
\emph{$d$-degenerate}\footnote{We remark that this matches the usual definition of $d$-degeneracy, i.e., that every subgraph has minimum degree at most $d$).}.
From the described interpretation, it is clear that if a graph is $d$-degenerate, then it is also $(d+1)$-list-colorable, and essentially the same argument also shows that it is $(d+1)$-correspondence
colorable.  Moreover, it is also easy to argue that it has Alon-Tarsi number at most $d+1$, and thus it is also online $(d+1)$-list-colorable.

However, degeneracy by itself is not sufficient for the results we discussed above, or even to show that planar graphs are 5-list-colorable.
Indeed, there exist planar graphs of minimum degree five, which are therefore not 4-degenerate.  However, it turns out that
it suffices to enhance the degeneracy ``just a bit'', and in fact we now know several different ways in which to do this; see~\cite{zhou2023arc} for an overview of some of the variations.
In this paper, we consider the notion of \emph{weak degeneracy}, which was proposed by Bernshteyn and Lee~\cite{weakdegen} and is one of the weakest of these notions,
in the sense that basically all other notions of generalized degeneracy are implied by weak degeneracy.  This has the advantage that the fact that a graph is weakly degenerate
has the widest range of consequences, but it also makes weak degeneracy the hardest to prove (for example, while Thomassen's 5-list-colorability proof can be
expressed in terms of several of the generalized degeneracy notions, it does not work in the weak degeneracy setting\footnote{%
We remark that Bernshteyn and Lee~\cite{weakdegen} claim otherwise in their original paper; however, as they later acknowledge~\cite{bernshteyn2024weak},
there is a subtle oversight in their argument.}).

To define weak degeneracy, let us consider the following slight refinement of the $\del$ operation.  We again consider a graph $G$ and a function $f:V(G)\to\mathbb{Z}$.
For an edge $vw$ of $G$, we define $\delsave(G, f, v, w) := (G-v, f'')$, where $f'':V(G-v) \rightarrow \mathbb{Z}$ is given by
$$f''(u) = 
   \begin{cases*}
   f(u)-1 & if $u\in N_G(v)\setminus \{w\}$ \\
   f(u) & otherwise.
\end{cases*}
$$
We say that this operation is \emph{legal} if $f(v)\ge 0$ and $f(v) > f(w)$.  Intuitively, this reflects the situation after we first choose a set $C\subseteq L(w)$ of exactly $f(w)+1$ colors
and a color $c\in L(v)\setminus C$ (which is possible when the operation is legal), use it to color $v$, and remove it from the lists of neighbors of $v$.
The choice of $C$ and $c$ ensures that the new list of $w$ contains all of $C$, and thus it still has length strictly larger than $f(w)$, justifying the correctness of defining $f''(w)$ as $f(w)$.
We say that $G$ is \emph{weakly $f$-degenerate} if there exists a sequence $(G,f)=(G_0,f_0)$, $(G_1,f_1)$, \ldots, $(G_n,f_n)$ such that for $i\in \{1,\ldots, n\}$,
either
\begin{itemize}
\item $(G_i,f_i)=\del(G_{i-1},f_{i-1},v_i)$ for a vertex $v_i\in V(G_{i-1})$, or
\item $(G_i,f_i)=\delsave(G_{i-1},f_{i-1},v_i,w_i)$ for an edge $v_iw_i\in E(G_{i-1})$,
\end{itemize}
and in either case the operation is legal.  If $f$ assigns the same value $d$ to all vertices of $G$, then we say that $G$ is \emph{weakly $d$-degenerate}.
It is again clear that if a graph is weakly $d$-degenerate, then it is also $(d+1)$-list-colorable, and even $(d+1)$-correspondence-colorable.  Moreover, it is also easy to argue that it has Alon-Tarsi number at most $d+1$ (see~\cite{zhou2023arc}), and thus it is also online $(d+1)$-list-colorable (which can be also shown directly~\cite{weakdegen}).

As we have mentioned above, Thomassen's elegant 5-list-colorability argument unfortunately does not seem to work in the weak degeneracy setting.
However, Dvořák, Lidický and Mohar~\cite{dvovrak20175} found a substantially more tedious (but less fragile) proof of 5-list-colorability of planar graphs,
and Bernshteyn, Lee, and Smith-Roberge~\cite{bernshteyn2024weak} were able to modify it to show that planar graphs are weakly $4$-degenerate.
This brings us to the question: Which other results can be recast in the weak degeneracy setting?

In particular, it is natural to ask whether proving
the existence of exponentially many colorings is possible via weak degeneracy.  This question was also raised in~\cite{davies2025localweakdegeneracyplanar}, where Davies and Smith-Roberge
propose the following natural formalism.  The \emph{availability} $a(O)$ of a legal operation $O$ is defined as
$$a(O)=\begin{cases*}
f(v)+1& if $O=\del(G, f, v)$\\
f(v)-f(w)& if $O=\delsave(G, f, v,w).$
\end{cases*}$$
That is, $a(O)$ gives a lower bound on the number of ways we can choose the color for $v$ if we follow the interpretation of the operations which we described above.
The \emph{total availability} of a sequence of operations $O_1$, \ldots, $O_n$ is the product $a=\prod_{i=1}^n a(O_i)$ of their availabilities,
and its \emph{average availability} is their geometric mean $\sqrt[n]{a}$.
We say that a graph is \emph{$\alpha$-availably weakly $f$-degenerate} (or \emph{totally $a$-availably weakly $f$-degenerate})
if there exists a sequence of legal operations as in the definition of weak $f$-degeneracy whose average availability is at least $\alpha$ (or whose total availability is at least $a$, respectively).
When $f$ assigns the same value $d$ to each vertex, we instead say that the graph is $\alpha$-availably or totally $a$-availably weakly $d$-degenerate.
Clearly, the following claim holds.
\begin{obs}
Let $G$ be a graph with $n$ vertices, let $f:V(G)\to\mathbb{Z}$ be a function, and let $\alpha\ge 1$ be a real number.
If $G$ is $\alpha$-availably weakly $f$-degenerate, then for every list assignment $L$ for $G$ such that $|L(v)|>f(v)$ for every $v\in V(G)$,
the number of $L$-colorings of $G$ is at least $\alpha^n$.
\end{obs}
The same observation also holds for correspondence coloring (the interpretation for online list coloring and for Alon-Tarsi number is less natural).
The main result of this paper now can be stated as follows.
\begin{thm}\label{thm:main}
Every planar graph is $2^{1/\CCC}$-availably weakly $4$-degenerate.
\end{thm}
In particular, this implies Theorem~\ref{thm:expmany}, even in the correspondence coloring setting.
Moreover, let us note the following consequence of $\alpha$-available weak $f$-degeneracy, which also holds in the correspondence coloring setting.
\begin{obs}\label{obs:remcol}
Let $\alpha\ge 1$ be a real number and $d\ge 2$ an integer.
If a graph $G$ is $\alpha$-availably weakly $d$-degenerate, then there exists a set $R\subseteq V(G)$ of size at least
$$(\log_{d+1}\alpha)\cdot |V(G)|$$
such that $G$ is $L$-colorable for every assignment $L$ of lists of size $d$ to vertices of $R$ and lists of size $d+1$ to the rest of the vertices of $G$.
\end{obs}
\begin{proof}
Let $O_1$, \ldots, $O_n$ be a sequence of operations that witnesses that $G$ is $\alpha$-availably weakly $d$-degenerate.
For $i\in\{1,\ldots,n\}$, let $v_i$ be the vertex removed in the operation $O_i$, and let
$$R=\{v_i:i\in\{1,\ldots,n\},a(O_i)>1\}.$$
Note that $a(O_i)\le d+1$ for each $i$, and thus
$$\alpha=\sqrt[n]{\prod_{i=1}^n a(O_i)}\le \sqrt[n]{(d+1)^{|R|}}$$
and
$$|R|\ge (\log_{d+1}\alpha)\cdot |V(G)|.$$
Consider any assignment $L$ of lists of size $d$ to vertices of $R$ and lists of size $d+1$ to the rest of the vertices of $G$.
Let $L'$ be the $(d+1)$-list-assignment obtained from $L$ by adding a new ``fake'' color to the list of each vertex of $R$.
We now proceed with $L'$-coloring the graph $G$ along the sequence $O_1$, \ldots, $O_n$ as we described at the definition
of $\del$ and $\delsave$ operations.  When we process an operation $O_i$ such that $v_i\in R$, we have at least $a(O_i)>1$
choices for the color we assign to $v_i$, and thus we can avoid using the fake color.  Therefore, the resulting coloring of $G$
will indeed be an $L$-coloring.
\end{proof}

Voigt~\cite{voigt1993list} proved that there exist planar graphs that are not $4$-list-colorable.  In contrast,
Theorem~\ref{thm:main} together with Observation~\ref{obs:remcol} shows that we can afford to have lists of size four on
at least $0.06n$ of its $n$ vertices.  Let us remark that this does not immediately (without a closer inspection) follow from the previous arguments
that show that planar graphs have exponentially many colorings from $5$-list-assignments, since in them the set of vertices where
the argument has a choice between several colors can depend on the contents of the lists, and thus it cannot be fixed in advance.


\subsection{Notation, terminology, and a technical strengthening of Theorem~\ref{thm:main}}

The following monotonicity properties are known to hold for weak degeneracy~\cite{weakdegen},
and essentially the same argument shows that they also hold in the availability setting.
\begin{obs}\label{obs:monot}
Let $G$ be a graph, let $f:V(G)\to \mathbb{Z}$ be a function, and let $a$ be a positive integer.
\begin{enumerate}[(i)]
\item\label{mon:adde} For all distinct non-adjacent vertices $u,v\in V(G)$, if the graph $G+uv$ is totally $a$-availably
weakly $f$-degenerate, then the graph $G$ is totally $a$-availably weakly $f$-degenerate as well.
\item\label{mon:decrf} Let $f':V(G)\to\mathbb{Z}$ be a function such that $f'(v)\le f(v)$ holds for every vertex $v\in V(G)$.
If the graph $G$ is totally $a$-availably weakly $f'$-degenerate, then it is also totally $a$-availably weakly $f$-degenerate.
\end{enumerate}
\end{obs}
\begin{proof}
We prove both claims by induction on the number of vertices.  In the first case, let $f'=f$;
and let $G'=G+uv$ in the first case and $G'=G$ in the second case.  By the definition of weak degeneracy with availability,
there exists an operation $O'$ ($\del(\ldots,x)$ or $\delsave(\ldots,x,y)$) whose
execution on $G'$ and $f'$ is legal, has availability $\alpha\ge 1$, and results in a pair $(H,g')$, where the graph
$H$ is totally $(a/\alpha)$-availably weakly $g'$-degenerate.  Let us now discuss the two claims separately.
\begin{enumerate}[(i)]
\item Recall that in this case $(G',f')=(G+uv,f)$.  If $x\not\in \{u,v\}$, then the execution of $O'$ on $(G,f)$ is legal, has availability $\alpha$, and results in
$(H-uv,g')$.  By the induction hypothesis (for $H-uv$, $H$, and $g'$, subcase \ref{mon:adde}), the graph $H-uv$ is totally $(a/\alpha)$-availably weakly $g'$-degenerate,
and thus $G$ is totally $a$-availably weakly $f$-degenerate.

Hence, suppose that $x\in \{u,v\}$, say $x=u$.  If $O'=\delsave(\ldots,u,v)$, then let $O=\del(\ldots,u)$, otherwise
let $O=O'$.  Since $uv\not\in E(G)$, the execution of $O$ on $G$ is legal, has availability at least $\alpha$,
and results in a pair $(H,g)$, where the function $g$ satisfies $g'(z)\le g(z)$ for every $z\in V(H)$.  By the induction hypothesis (for $H$, $g$, and $g'$,
subcase \ref{mon:decrf}), the graph $H$ is totally $(a/\alpha)$-availably weakly $g$-degenerate,
and thus $G$ is totally $a$-availably weakly $f$-degenerate.
\item Recall that in this case $G'=G$ and $f'(v)\le f(v)$ holds for every vertex $v\in V(G)$.  If $O'$ is $\delsave(\ldots,x,y)$ and $f(y)>f'(y)$, then let $O=\del(\ldots,x)$, otherwise let $O=O'$.
Since $f(x)\ge f'(x)$, this ensures that the operation $O$ is legal for $(G,f)$ and has availability at least $\alpha$.
Moreover, the result of the operation is a pair $(H,g)$, where the function $g$ satisfies $g'(z)\le g(z)$ for every $z\in V(H)$---in
the case where $O\neq O'$, we have $g'(y)=f'(y)\le f(y)-1=g(y)$, since $f(y)>f'(y)$.
By the induction hypothesis (for $H$, $g$, and $g'$, subcase \ref{mon:decrf}),
the graph $H$ is totally $(a/\alpha)$-availably weakly $g$-degenerate,
and thus $G$ is totally $a$-availably weakly $f$-degenerate.
\end{enumerate}
\end{proof}

It will be convenient to add one more (redundant) operation to the definition of weak degeneracy.
For a graph $G$, a function $f:V(G)\to\mathbb{Z}$, and a vertex $v\in V(G)$,
we let $\Omit(G,f,v)=(G-v,f')$, where $f'$ is the restriction of $f$ to $V(G)\setminus\{v\}$.
We say that this operation is \emph{legal} if $f(v)\ge \deg_G v$, and the \emph{availability} of this
operation is defined as $f(v)-\deg_G v+1$.  Note that any sequence of legal $\del$, $\delsave$, and $\Omit$
operations that removes the whole graph can be transformed into a sequence of legal $\del$ and $\delsave$ operations with the same effect and
availability: To eliminate the last $\Omit(G,f,v)$ operation from the sequence, we ignore it and perform the remaining operations, obtaining
a pair $(G_0,f_0)$, where $V(G_0)=\{v\}$.  Since all neighbors of $v$ were removed by $\del$'s and $\delsave$'s, we have $f_0(v)=f(v)-\deg_G v\ge 0$.  We then add the operation
$\del(G_0,f_0,v)$, whose availability is exactly $f_0(v)+1=f(v)-\deg_G v+1$, i.e., equal to the availability of the eliminated $\Omit$ operation.
By repeating this procedure, we eventually eliminate all $\Omit$ operations from the sequence, without affecting its total availability.

Thus, we can use $\Omit$ operations when convenient, without affecting the notions of weak degeneracy or availability.
Moreover, by the described transformation, we can always assume that a sequence of operations showing that a graph is weakly degenerate only uses
the $\del$ and $\delsave$ operations (to emphasize this, we call such a sequence \emph{$\Omit$-free}).

Let $v(G)$ denote the number of vertices of a graph $G$.  The \emph{length} of a path $P$ is defined as its number of edges,
i.e., $v(P)-1$.
Given a plane graph $G$, we refer to the infinite face of $G$ as its
\emph{outer face}, and denote by $\partial G$ the graph whose vertex- and
edge-sets are precisely those in the boundary of the outer face of $G$.
All other faces of $G$ are \emph{internal}.
The \emph{length} of a face $g$ is the sum of the lengths of its boundary walks (and in particular, the bridges
contained in the boundary of $g$ contribute $2$ to the length of $g$).

A path $P$ in a plane graph $G$ is a \emph{generalized chord} if its endpoints are in $V(\partial G)$ and $P$ is otherwise disjoint from $\partial G$
(and in particular $E(P)\cap E(\partial G)=\emptyset$), and in case that $P$ consists of a single vertex $p$, then additionally $p$
is a cutvertex of $\partial G$.  Note that if $G$ is a generalized chord, then there exists a simple closed curve
$\gamma$ intersecting the drawing of $G$ exactly along $P$ and otherwise drawn in the outer face of $G$ such that at least one
vertex or edge of $G$ is drawn in each of the two open regions $\rho_1$ and $\rho_2$ to which $\gamma$ divides the plane.
For $i\in \{1,2\}$, let $G_i$ be the subgraph of $G$ drawn in the closure of $\rho_i$, where we view $G_i$ as a plane
graph with the drawing inherited from $G$; we say that $P$ \emph{separates $G$ into $G_1$ and $G_2$}.  Note that $G_1 \cap G_2 = P =\partial G_1\cap \partial G_2$, $\partial G_1\neq P\neq\partial G_2$,
$G_1 \cup G_2 = G$, and $\partial G_1\cup \partial G_2=P\cup \partial G$.
Let us remark that if $\partial G$ is a cycle, then $G_1$ and $G_2$
are determined by $P$ uniquely up to their exchange, but otherwise (when one of the ends of $P$ is a cutvertex of $\partial G$),
there may be several possibilities how to choose $G_1$ and $G_2$.  When a generalized chord $P$ has length $t$, we say that $P$ is
a \emph{$t$-chord}.

Following the standard approach in the area, instead of proving Theorem~\ref{thm:main} directly, we prove
its strengthening with a tighter degeneracy condition on the vertices incident with the outer face.
Thus, the basic object we work with is as follows.

\begin{definition}
A \emph{canvas} is a triple $K=(G,P,f)$, where
\begin{itemize}
    \item $G$ is a plane graph,
    \item $P$ is a path of length at most two in the outer face boundary walk of $G$, and
    \item $f$ is an integer-valued function whose domain includes $V(G)\setminus V(P)$.
\end{itemize}
The canvas is \emph{valid} if $P$ is an induced path in $G$ and
    \begin{itemize}
    \item $f(v) \geq 4$ for all $v \in V(G)\setminus V(\partial G)$,
    \item $f(v) \geq 2$ for all $v \in V(\partial G) \setminus V(P)$, and
    \item $f^{-1}(2)\cap V(\partial G) \setminus V(P)$ is an independent set in $G$.
    \end{itemize}
We write $f_K$ for the function $f_K:V(G)\setminus V(P)\to \mathbb Z$ given by $f_K(v) = f(v) - |N(v)\cap V(P)|$
for each $v\in V(G)\setminus V(P)$, and we say that the canvas $K$ is \emph{($\alpha$-availably, resp. totally $a$-availably) weakly degenerate} if $P$ is an induced path in $G$ and
the graph $G-V(P)$ is ($\alpha$-availably, resp. totally $a$-availably) weakly $f_K$-degenerate.
\end{definition}

As an intuition for the definitions, suppose that $G_0$ is a plane graph and $f:V(G_0)\to\mathbb{Z}$
is a function such that $f(v)\ge 4$ for all vertices $v\in V(G_0)$ not incident with the outer face, $f(v)\ge 2$ for
all $v\in V(\partial G_0)$, and $f^{-1}(2)$ is a subset of $V(\partial G)\setminus V(P)$ that is an independent set in $G_0$.  Let $P$ be a generalized chord in $G_0$ of length at most
two separating $G_0$ into $H$ and $G$ such that $P$ is an induced path in $G$.  Then $K=(G,P,f)$ is a valid canvas,
and moreover if the graph $H$ is ($\alpha$-availably) weakly $f$-degenerate and the canvas $K$ is ($\alpha$-availably) weakly degenerate, then the whole graph $G_0$
is ($\alpha$-availably) weakly $f$-degenerate.  Indeed, we can first perform the sequence of operations showing that $H$ is ($\alpha$-availably) weakly $f$-degenerate on $G_0$,
resulting in the pair $(G-V(P),f_K)$, then perform the sequence of operations showing that $K$ is ($\alpha$-availably) weakly degenerate.

We lift the operations $\del$, $\delsave$, and $\Omit$ to canvases: Consider a canvas $K=(G,P,f)$ and vertices $u,v\in V(G)\setminus V(P)$.
We let $\del(K,v)$ be the canvas $(G',P,f')$ where $(G',f')=\del(G,f,v)$, let $\delsave(K,v,u)$ be the canvas $(G'',P,f'')$
where $(G'',f'')=\delsave(G,f,v,u)$, and let $\Omit(K,v)$ be the canvas $(G''',P,f''')$ where $(G''',f''')=\Omit(G,f,v)$
(that is, we decrease the values of $f$ exactly as when performing these operations on graphs).
We say that these operations are \emph{legal} for $K$ if they are legal in the graph $G-V(P)$
using the function $f_K$, and we define their \emph{availability} as the
availability of the corresponding operations for $G-V(P)$ and $f_K$.
For example, the operation $\del(K,v)$ is legal if $f(v)\ge |N(v)\cap V(P)|$, and its availability is $f(v)-|N(v)\cap V(P)|+1$.
Let us remark that we do not require the resulting canvases to be valid (even if $K$ itself is valid);
indeed, we often perform sequences of operations where the initial and the final canvases are valid, but the intermediate
ones are not.

Not all valid canvases are weakly degenerate.  The \emph{tripod} is the canvas $W=(H,P,f)$
such that $P$ is an induced path of length two in $H$, $f$ has domain $V(H)\setminus V(P)$, $V(H) \setminus V(P)$ consists of a single vertex $v$ with $f(v) = 2$,
and $v$ is adjacent to all three vertices of $P$.  Clearly a tripod is valid, but not weakly degenerate, since $f_W(v)=-1$.
Moreover, canvases that contain a tripod as a subcanvas are not weakly degenerate (a canvas $K' = (G',P',f')$
is a \emph{subcanvas} of a canvas $K = (G,P,f)$ if $G'\subseteq G$, $P'=P$, and $f'$ and $f$ have the same values
on the vertices of $V(G')\setminus V(P')$).  A canvas is \emph{tripod-free} if it does not contain the tripod as a subcanvas.

\begin{figure}[ht]\centering
    \subcaptionbox{\ref{S}-singular}[0.3\textwidth]{
\begin{tikzpicture}[scale=1]
    \node[precoloured] (u1) at (120:2cm){};
    \node[precoloured] (u2) at (90:2cm){};
    \node[precoloured] (u3) at (60:2cm){};
    \node[4vertex] (v) at (0,0) {};
    \draw (v) -- (u1) -- (u2) -- (u3) -- (v) -- (u2);
\end{tikzpicture}


\begin{tikzpicture}[scale=1]
    \node[precoloured] (u1) at (120:2cm){};
    \node[precoloured] (u2) at (90:2cm){};
    \node[precoloured] (u3) at (60:2cm){};
    \node[3vertex] (v) at (0,0) {};
    \node[geq4vertex] (u) at (0,1) {};
    \draw (v) -- (u1) -- (u2) -- (u3) -- (v);
    \draw (u2) -- (u) -- (u1);
    \draw (v) -- (u) -- (u3);
\end{tikzpicture}

\begin{tikzpicture}[scale=1]
    \node[precoloured] (u1) at (120:2cm){};
    \node[precoloured] (u2) at (90:2cm){};
    \node[precoloured] (u3) at (60:2cm){};
    \node[4vertex] (v) at (0.5,0) {};
    \node[3vertex] (u) at (-0.5,0) {};

    \draw (u) -- (u1) -- (u2) -- (u3);
    \draw (v) -- (u) -- (u2) -- (v);
    \draw[wavy] (u3) -- (v);
\end{tikzpicture}

\begin{tikzpicture}[scale=1]
    \node[precoloured] (u1) at (120:2cm){};
    \node[precoloured] (u2) at (90:2cm){};
    \node[precoloured] (u3) at (60:2cm){};
    \node[3vertex] (v1) at (-0.6,0) {};
    \node[4vertex] (v2) at (0,0) {};
    \node[3vertex] (v3) at (0.6,0) {};
    \draw (v1) -- (u1) -- (u2) -- (u3) -- (v3) -- (v2) -- (v1);
    \draw (v1) -- (u2) -- (v2);
    \draw (u2) -- (v3);
\end{tikzpicture}

\begin{tikzpicture}[scale=1]
    \node[precoloured] (u1) at (120:2cm){};
    \node[precoloured] (u2) at (60:2cm){};
    \node[3vertex] (v) at (0,0) {};
    \draw (v) -- (u1) -- (u2) -- (v);
\end{tikzpicture}

\begin{tikzpicture}[scale=1]
    \node[precoloured] (u1) at (120:2cm){};
    \node[precoloured] (u2) at (90:2cm){};
    \node[precoloured] (u3) at (60:2cm){};
    \node[4vertex] (v) at (0,0) {};
    \node[3vertex] (u) at (1,0.2) {};
    \draw (v) -- (u1) -- (u2) -- (u3) -- (v) -- (u2);
    \draw (u3) -- (u) -- (v);
\end{tikzpicture}

\begin{tikzpicture}[scale=1]
    \node[precoloured] (u1) at (120:2cm){};
    \node[precoloured] (u2) at (90:2cm){};
    \node[precoloured] (u3) at (60:2cm){};
    \node[4vertex] (v) at (0.5,0) {};
    \node[3vertex] (u) at (-0.5,0) {};
    \node[3vertex] (w) at (1.5,0.2) {};

    \draw (u) -- (u1) -- (u2) -- (u3) -- (v) -- (u) -- (u2) -- (v);
    \draw (u3) -- (w) -- (v);
\end{tikzpicture}
\begin{tikzpicture}[scale=1]
    \node[precoloured] (u1) at (120:2cm){};
    \node[precoloured] (u2) at (90:2cm){};
    \node[precoloured] (u3) at (60:2cm){};
    \node[4vertex] (v) at (0,0) {};
    \node[3vertex] (u) at (1,0.2) {};
    \node[3vertex] (w) at (-1,0.2) {};

    \draw (v) -- (u1) -- (u2) -- (u3) -- (v) -- (u2);
    \draw (u3) -- (u) -- (v);
    \draw (u1) -- (w) -- (v);

\end{tikzpicture}}%
    \hspace*{0.01\textwidth}
    \subcaptionbox{\ref{S3}-singular}[0.3\textwidth]{
\begin{tikzpicture}[scale=1]
    \node[precoloured] (u1) at (120:2cm){};
    \node[precoloured] (u2) at (90:2cm){};
    \node[precoloured] (u3) at (60:2cm){};
    \node[4vertex] (v) at (0,0) {};
    \draw (v) -- (u1) -- (u2) -- (u3) -- (v) -- (u2);
\end{tikzpicture}


\begin{tikzpicture}[scale=1]
    \node[precoloured] (u1) at (120:2cm){};
    \node[precoloured] (u2) at (90:2cm){};
    \node[precoloured] (u3) at (60:2cm){};
    \node[3vertex] (v) at (0,0) {};
    \node[geq4vertex] (u) at (0,1) {};
    \draw (v) -- (u1) -- (u2) -- (u3) -- (v);
    \draw (u2) -- (u) -- (u1);
    \draw (v) -- (u) -- (u3);
\end{tikzpicture}

\begin{tikzpicture}[scale=1]
    \node[precoloured] (u1) at (120:2cm){};
    \node[precoloured] (u2) at (90:2cm){};
    \node[precoloured] (u3) at (60:2cm){};
    \node[4vertex] (v) at (0.5,0) {};
    \node[3vertex] (u) at (-0.5,0) {};

    \draw (u) -- (u1) -- (u2) -- (u3);
    \draw (v) -- (u) -- (u2) -- (v);
    \draw[wavy] (u3) -- (v);
\end{tikzpicture}

\begin{tikzpicture}[scale=1]
    \node[precoloured] (u1) at (120:2cm){};
    \node[precoloured] (u2) at (90:2cm){};
    \node[precoloured] (u3) at (60:2cm){};
    \node[3vertex] (v1) at (-0.6,0) {};
    \node[4vertex] (v2) at (0,0) {};
    \node[3vertex] (v3) at (0.6,0) {};
    \draw (v1) -- (u1) -- (u2) -- (u3) -- (v3) -- (v2) -- (v1);
    \draw (v1) -- (u2) -- (v2);
    \draw (u2) -- (v3);
\end{tikzpicture}

\begin{tikzpicture}[scale=1]
    \node[precoloured] (u1) at (120:2cm){};
    \node[precoloured] (u2) at (60:2cm){};
    \node[3vertex] (v) at (0,0) {};
    \draw (v) -- (u1) -- (u2) -- (v);
\end{tikzpicture}

\begin{tikzpicture}[scale=1]
    \node[precoloured] (u1) at (120:2cm){};
    \node[precoloured] (u2) at (90:2cm){};
    \node[precoloured] (u3) at (60:2cm){};
    \node[4vertex] (v) at (0,0) {};
    \node[3vertex] (u) at (1,0.2) {};
    \draw (v) -- (u1) -- (u2) -- (u3) -- (v) -- (u2);
    \draw (u3) -- (u) -- (v);
\end{tikzpicture}

\begin{tikzpicture}[scale=1]
    \node[precoloured] (u1) at (120:2cm){};
    \node[precoloured] (u2) at (90:2cm){};
    \node[precoloured] (u3) at (60:2cm){};
    \node[4vertex] (v) at (0.5,0) {};
    \node[3vertex] (u) at (-0.5,0) {};
    \node[3vertex] (w) at (1.5,0.2) {};

    \draw (u) -- (u1) -- (u2) -- (u3) -- (v) -- (u) -- (u2) -- (v);
    \draw (u3) -- (w) -- (v);
\end{tikzpicture}
\begin{tikzpicture}[scale=1]
    \node[precoloured] (u1) at (120:2cm){};
    \node[precoloured] (u2) at (90:2cm){};
    \node[precoloured] (u3) at (60:2cm){};
    \node[4vertex] (v) at (0,0) {};
    \node[3vertex] (u) at (1,0.2) {};
    \node[3vertex] (w) at (-1,0.2) {};

    \draw (v) -- (u1) -- (u2) -- (u3) -- (v) -- (u2);
    \draw (u3) -- (u) -- (v);
    \draw (u1) -- (w) -- (v);

\end{tikzpicture}}%
    \hspace*{0.01\textwidth}
    \subcaptionbox{\ref{S2i}-singular}[0.3\textwidth]{
\begin{tikzpicture}[scale=1]
    \node[precoloured] (u1) at (120:2cm){};
    \node[precoloured] (u2) at (90:2cm){};
    \node[precoloured] (u3) at (60:2cm){};
    \node[4vertex] (v) at (0,0) {};
    \draw (v) -- (u1) -- (u2) -- (u3) -- (v) -- (u2);
\end{tikzpicture}


\begin{tikzpicture}[scale=1]
    \node[precoloured] (u1) at (120:2cm){};
    \node[precoloured] (u2) at (90:2cm){};
    \node[precoloured] (u3) at (60:2cm){};
    \node[3vertex] (v) at (0,0) {};
    \node[geq4vertex] (u) at (0,1) {};
    \draw (v) -- (u1) -- (u2) -- (u3) -- (v);
    \draw (u2) -- (u) -- (u1);
    \draw (v) -- (u) -- (u3);
\end{tikzpicture}

\begin{tikzpicture}[scale=1]
    \node[precoloured] (u1) at (120:2cm){};
    \node[precoloured] (u2) at (90:2cm){};
    \node[precoloured] (u3) at (60:2cm){};
    \node[4vertex] (v) at (0.5,0) {};
    \node[3vertex] (u) at (-0.5,0) {};

    \draw (u) -- (u1) -- (u2) -- (u3);
    \draw (v) -- (u) -- (u2) -- (v);
    \draw[wavy] (u3) -- (v);
\end{tikzpicture}

\begin{tikzpicture}[scale=1]
    \node[precoloured] (u1) at (120:2cm){};
    \node[precoloured] (u2) at (90:2cm){};
    \node[precoloured] (u3) at (60:2cm){};
    \node[3vertex] (v1) at (-0.6,0) {};
    \node[4vertex] (v2) at (0,0) {};
    \node[3vertex] (v3) at (0.6,0) {};
    \draw (v1) -- (u1) -- (u2) -- (u3) -- (v3) -- (v2) -- (v1);
    \draw (v1) -- (u2) -- (v2);
    \draw (u2) -- (v3);
\end{tikzpicture}

\begin{tikzpicture}[scale=1]
    \node[precoloured] (u1) at (120:2cm){};
    \node[precoloured] (u2) at (60:2cm){};
    \node[3vertex] (v) at (0,0) {};
    \draw (v) -- (u1) -- (u2) -- (v);
\end{tikzpicture}

\begin{tikzpicture}[scale=1]
    \node[precoloured] (u1) at (120:2cm){};
    \node[precoloured] (u2) at (90:2cm){};
    \node[precoloured] (u3) at (60:2cm){};
    \node[4vertex] (v) at (0,0) {};
    \node[3vertex] (u) at (1,0.2) {};
    \draw (v) -- (u1) -- (u2) -- (u3) -- (v) -- (u2);
    \draw (u3) -- (u) -- (v);
\end{tikzpicture}

\begin{tikzpicture}[scale=1]
    \node[precoloured] (u1) at (120:2cm){};
    \node[precoloured] (u2) at (90:2cm){};
    \node[precoloured] (u3) at (60:2cm){};
    \node[4vertex] (v) at (0.5,0) {};
    \node[3vertex] (u) at (-0.5,0) {};
    \node[3vertex] (w) at (1.5,0.2) {};

    \draw (u) -- (u1) -- (u2) -- (u3) -- (v) -- (u) -- (u2) -- (v);
    \draw (u3) -- (w) -- (v);
\end{tikzpicture}
\begin{tikzpicture}[scale=1]
    \node[precoloured] (u1) at (120:2cm){};
    \node[precoloured] (u2) at (90:2cm){};
    \node[precoloured] (u3) at (60:2cm){};
    \node[4vertex] (v) at (0,0) {};
    \node[3vertex] (u) at (1,0.2) {};
    \node[3vertex] (w) at (-1,0.2) {};

    \draw (v) -- (u1) -- (u2) -- (u3) -- (v) -- (u2);
    \draw (u3) -- (u) -- (v);
    \draw (u1) -- (w) -- (v);

\end{tikzpicture}}%

\vspace*{8mm}

    \subcaptionbox{\ref{S23a}-singular}[0.3\textwidth]{
\begin{tikzpicture}[scale=1]
    \node[precoloured] (u1) at (120:2cm){};
    \node[precoloured] (u2) at (90:2cm){};
    \node[precoloured] (u3) at (60:2cm){};
    \node[4vertex] (v) at (0,0) {};
    \draw (v) -- (u1) -- (u2) -- (u3) -- (v) -- (u2);
\end{tikzpicture}


\begin{tikzpicture}[scale=1]
    \node[precoloured] (u1) at (120:2cm){};
    \node[precoloured] (u2) at (90:2cm){};
    \node[precoloured] (u3) at (60:2cm){};
    \node[3vertex] (v) at (0,0) {};
    \node[geq4vertex] (u) at (0,1) {};
    \draw (v) -- (u1) -- (u2) -- (u3) -- (v);
    \draw (u2) -- (u) -- (u1);
    \draw (v) -- (u) -- (u3);
\end{tikzpicture}

\begin{tikzpicture}[scale=1]
    \node[precoloured] (u1) at (120:2cm){};
    \node[precoloured] (u2) at (90:2cm){};
    \node[precoloured] (u3) at (60:2cm){};
    \node[4vertex] (v) at (0.5,0) {};
    \node[3vertex] (u) at (-0.5,0) {};

    \draw (u) -- (u1) -- (u2) -- (u3);
    \draw (v) -- (u) -- (u2) -- (v);
    \draw[wavy] (u3) -- (v);
\end{tikzpicture}

\begin{tikzpicture}[scale=1]
    \node[precoloured] (u1) at (120:2cm){};
    \node[precoloured] (u2) at (90:2cm){};
    \node[precoloured] (u3) at (60:2cm){};
    \node[3vertex] (v1) at (-0.6,0) {};
    \node[4vertex] (v2) at (0,0) {};
    \node[3vertex] (v3) at (0.6,0) {};
    \draw (v1) -- (u1) -- (u2) -- (u3) -- (v3) -- (v2) -- (v1);
    \draw (v1) -- (u2) -- (v2);
    \draw (u2) -- (v3);
\end{tikzpicture}

\begin{tikzpicture}[scale=1]
    \node[precoloured] (u1) at (120:2cm){};
    \node[precoloured] (u2) at (60:2cm){};
    \node[3vertex] (v) at (0,0) {};
    \draw (v) -- (u1) -- (u2) -- (v);
\end{tikzpicture}

\begin{tikzpicture}[scale=1]
    \node[precoloured] (u1) at (120:2cm){};
    \node[precoloured] (u2) at (90:2cm){};
    \node[precoloured] (u3) at (60:2cm){};
    \node[4vertex] (v) at (0,0) {};
    \node[3vertex] (u) at (1,0.2) {};
    \draw (v) -- (u1) -- (u2) -- (u3) -- (v) -- (u2);
    \draw (u3) -- (u) -- (v);
\end{tikzpicture}

\begin{tikzpicture}[scale=1]
    \node[precoloured] (u1) at (120:2cm){};
    \node[precoloured] (u2) at (90:2cm){};
    \node[precoloured] (u3) at (60:2cm){};
    \node[4vertex] (v) at (0.5,0) {};
    \node[3vertex] (u) at (-0.5,0) {};
    \node[3vertex] (w) at (1.5,0.2) {};

    \draw (u) -- (u1) -- (u2) -- (u3) -- (v) -- (u) -- (u2) -- (v);
    \draw (u3) -- (w) -- (v);
\end{tikzpicture}
\begin{tikzpicture}[scale=1]
    \node[precoloured] (u1) at (120:2cm){};
    \node[precoloured] (u2) at (90:2cm){};
    \node[precoloured] (u3) at (60:2cm){};
    \node[4vertex] (v) at (0,0) {};
    \node[3vertex] (u) at (1,0.2) {};
    \node[3vertex] (w) at (-1,0.2) {};

    \draw (v) -- (u1) -- (u2) -- (u3) -- (v) -- (u2);
    \draw (u3) -- (u) -- (v);
    \draw (u1) -- (w) -- (v);

\end{tikzpicture}}%
    \subcaptionbox{\ref{S232a}-singular}[0.3\textwidth]{
\begin{tikzpicture}[scale=1]
    \node[precoloured] (u1) at (120:2cm){};
    \node[precoloured] (u2) at (90:2cm){};
    \node[precoloured] (u3) at (60:2cm){};
    \node[4vertex] (v) at (0,0) {};
    \draw (v) -- (u1) -- (u2) -- (u3) -- (v) -- (u2);
\end{tikzpicture}


\begin{tikzpicture}[scale=1]
    \node[precoloured] (u1) at (120:2cm){};
    \node[precoloured] (u2) at (90:2cm){};
    \node[precoloured] (u3) at (60:2cm){};
    \node[3vertex] (v) at (0,0) {};
    \node[geq4vertex] (u) at (0,1) {};
    \draw (v) -- (u1) -- (u2) -- (u3) -- (v);
    \draw (u2) -- (u) -- (u1);
    \draw (v) -- (u) -- (u3);
\end{tikzpicture}

\begin{tikzpicture}[scale=1]
    \node[precoloured] (u1) at (120:2cm){};
    \node[precoloured] (u2) at (90:2cm){};
    \node[precoloured] (u3) at (60:2cm){};
    \node[4vertex] (v) at (0.5,0) {};
    \node[3vertex] (u) at (-0.5,0) {};

    \draw (u) -- (u1) -- (u2) -- (u3);
    \draw (v) -- (u) -- (u2) -- (v);
    \draw[wavy] (u3) -- (v);
\end{tikzpicture}

\begin{tikzpicture}[scale=1]
    \node[precoloured] (u1) at (120:2cm){};
    \node[precoloured] (u2) at (90:2cm){};
    \node[precoloured] (u3) at (60:2cm){};
    \node[3vertex] (v1) at (-0.6,0) {};
    \node[4vertex] (v2) at (0,0) {};
    \node[3vertex] (v3) at (0.6,0) {};
    \draw (v1) -- (u1) -- (u2) -- (u3) -- (v3) -- (v2) -- (v1);
    \draw (v1) -- (u2) -- (v2);
    \draw (u2) -- (v3);
\end{tikzpicture}

\begin{tikzpicture}[scale=1]
    \node[precoloured] (u1) at (120:2cm){};
    \node[precoloured] (u2) at (60:2cm){};
    \node[3vertex] (v) at (0,0) {};
    \draw (v) -- (u1) -- (u2) -- (v);
\end{tikzpicture}

\begin{tikzpicture}[scale=1]
    \node[precoloured] (u1) at (120:2cm){};
    \node[precoloured] (u2) at (90:2cm){};
    \node[precoloured] (u3) at (60:2cm){};
    \node[4vertex] (v) at (0,0) {};
    \node[3vertex] (u) at (1,0.2) {};
    \draw (v) -- (u1) -- (u2) -- (u3) -- (v) -- (u2);
    \draw (u3) -- (u) -- (v);
\end{tikzpicture}

\begin{tikzpicture}[scale=1]
    \node[precoloured] (u1) at (120:2cm){};
    \node[precoloured] (u2) at (90:2cm){};
    \node[precoloured] (u3) at (60:2cm){};
    \node[4vertex] (v) at (0.5,0) {};
    \node[3vertex] (u) at (-0.5,0) {};
    \node[3vertex] (w) at (1.5,0.2) {};

    \draw (u) -- (u1) -- (u2) -- (u3) -- (v) -- (u) -- (u2) -- (v);
    \draw (u3) -- (w) -- (v);
\end{tikzpicture}
\begin{tikzpicture}[scale=1]
    \node[precoloured] (u1) at (120:2cm){};
    \node[precoloured] (u2) at (90:2cm){};
    \node[precoloured] (u3) at (60:2cm){};
    \node[4vertex] (v) at (0,0) {};
    \node[3vertex] (u) at (1,0.2) {};
    \node[3vertex] (w) at (-1,0.2) {};

    \draw (v) -- (u1) -- (u2) -- (u3) -- (v) -- (u2);
    \draw (u3) -- (u) -- (v);
    \draw (u1) -- (w) -- (v);

\end{tikzpicture}}%
    \subcaptionbox{\ref{S23b}-singular}[0.3\textwidth]{
\begin{tikzpicture}[scale=1]
    \node[precoloured] (u1) at (120:2cm){};
    \node[precoloured] (u2) at (90:2cm){};
    \node[precoloured] (u3) at (60:2cm){};
    \node[4vertex] (v) at (0,0) {};
    \draw (v) -- (u1) -- (u2) -- (u3) -- (v) -- (u2);
\end{tikzpicture}


\begin{tikzpicture}[scale=1]
    \node[precoloured] (u1) at (120:2cm){};
    \node[precoloured] (u2) at (90:2cm){};
    \node[precoloured] (u3) at (60:2cm){};
    \node[3vertex] (v) at (0,0) {};
    \node[geq4vertex] (u) at (0,1) {};
    \draw (v) -- (u1) -- (u2) -- (u3) -- (v);
    \draw (u2) -- (u) -- (u1);
    \draw (v) -- (u) -- (u3);
\end{tikzpicture}

\begin{tikzpicture}[scale=1]
    \node[precoloured] (u1) at (120:2cm){};
    \node[precoloured] (u2) at (90:2cm){};
    \node[precoloured] (u3) at (60:2cm){};
    \node[4vertex] (v) at (0.5,0) {};
    \node[3vertex] (u) at (-0.5,0) {};

    \draw (u) -- (u1) -- (u2) -- (u3);
    \draw (v) -- (u) -- (u2) -- (v);
    \draw[wavy] (u3) -- (v);
\end{tikzpicture}

\begin{tikzpicture}[scale=1]
    \node[precoloured] (u1) at (120:2cm){};
    \node[precoloured] (u2) at (90:2cm){};
    \node[precoloured] (u3) at (60:2cm){};
    \node[3vertex] (v1) at (-0.6,0) {};
    \node[4vertex] (v2) at (0,0) {};
    \node[3vertex] (v3) at (0.6,0) {};
    \draw (v1) -- (u1) -- (u2) -- (u3) -- (v3) -- (v2) -- (v1);
    \draw (v1) -- (u2) -- (v2);
    \draw (u2) -- (v3);
\end{tikzpicture}

\begin{tikzpicture}[scale=1]
    \node[precoloured] (u1) at (120:2cm){};
    \node[precoloured] (u2) at (60:2cm){};
    \node[3vertex] (v) at (0,0) {};
    \draw (v) -- (u1) -- (u2) -- (v);
\end{tikzpicture}

\begin{tikzpicture}[scale=1]
    \node[precoloured] (u1) at (120:2cm){};
    \node[precoloured] (u2) at (90:2cm){};
    \node[precoloured] (u3) at (60:2cm){};
    \node[4vertex] (v) at (0,0) {};
    \node[3vertex] (u) at (1,0.2) {};
    \draw (v) -- (u1) -- (u2) -- (u3) -- (v) -- (u2);
    \draw (u3) -- (u) -- (v);
\end{tikzpicture}

\begin{tikzpicture}[scale=1]
    \node[precoloured] (u1) at (120:2cm){};
    \node[precoloured] (u2) at (90:2cm){};
    \node[precoloured] (u3) at (60:2cm){};
    \node[4vertex] (v) at (0.5,0) {};
    \node[3vertex] (u) at (-0.5,0) {};
    \node[3vertex] (w) at (1.5,0.2) {};

    \draw (u) -- (u1) -- (u2) -- (u3) -- (v) -- (u) -- (u2) -- (v);
    \draw (u3) -- (w) -- (v);
\end{tikzpicture}
\begin{tikzpicture}[scale=1]
    \node[precoloured] (u1) at (120:2cm){};
    \node[precoloured] (u2) at (90:2cm){};
    \node[precoloured] (u3) at (60:2cm){};
    \node[4vertex] (v) at (0,0) {};
    \node[3vertex] (u) at (1,0.2) {};
    \node[3vertex] (w) at (-1,0.2) {};

    \draw (v) -- (u1) -- (u2) -- (u3) -- (v) -- (u2);
    \draw (u3) -- (u) -- (v);
    \draw (u1) -- (w) -- (v);

\end{tikzpicture}}%

\vspace*{5mm}

    \subcaptionbox{\ref{S232b}-singular}[0.4\textwidth]{
\begin{tikzpicture}[scale=1]
    \node[precoloured] (u1) at (120:2cm){};
    \node[precoloured] (u2) at (90:2cm){};
    \node[precoloured] (u3) at (60:2cm){};
    \node[4vertex] (v) at (0,0) {};
    \draw (v) -- (u1) -- (u2) -- (u3) -- (v) -- (u2);
\end{tikzpicture}


\begin{tikzpicture}[scale=1]
    \node[precoloured] (u1) at (120:2cm){};
    \node[precoloured] (u2) at (90:2cm){};
    \node[precoloured] (u3) at (60:2cm){};
    \node[3vertex] (v) at (0,0) {};
    \node[geq4vertex] (u) at (0,1) {};
    \draw (v) -- (u1) -- (u2) -- (u3) -- (v);
    \draw (u2) -- (u) -- (u1);
    \draw (v) -- (u) -- (u3);
\end{tikzpicture}

\begin{tikzpicture}[scale=1]
    \node[precoloured] (u1) at (120:2cm){};
    \node[precoloured] (u2) at (90:2cm){};
    \node[precoloured] (u3) at (60:2cm){};
    \node[4vertex] (v) at (0.5,0) {};
    \node[3vertex] (u) at (-0.5,0) {};

    \draw (u) -- (u1) -- (u2) -- (u3);
    \draw (v) -- (u) -- (u2) -- (v);
    \draw[wavy] (u3) -- (v);
\end{tikzpicture}

\begin{tikzpicture}[scale=1]
    \node[precoloured] (u1) at (120:2cm){};
    \node[precoloured] (u2) at (90:2cm){};
    \node[precoloured] (u3) at (60:2cm){};
    \node[3vertex] (v1) at (-0.6,0) {};
    \node[4vertex] (v2) at (0,0) {};
    \node[3vertex] (v3) at (0.6,0) {};
    \draw (v1) -- (u1) -- (u2) -- (u3) -- (v3) -- (v2) -- (v1);
    \draw (v1) -- (u2) -- (v2);
    \draw (u2) -- (v3);
\end{tikzpicture}

\begin{tikzpicture}[scale=1]
    \node[precoloured] (u1) at (120:2cm){};
    \node[precoloured] (u2) at (60:2cm){};
    \node[3vertex] (v) at (0,0) {};
    \draw (v) -- (u1) -- (u2) -- (v);
\end{tikzpicture}

\begin{tikzpicture}[scale=1]
    \node[precoloured] (u1) at (120:2cm){};
    \node[precoloured] (u2) at (90:2cm){};
    \node[precoloured] (u3) at (60:2cm){};
    \node[4vertex] (v) at (0,0) {};
    \node[3vertex] (u) at (1,0.2) {};
    \draw (v) -- (u1) -- (u2) -- (u3) -- (v) -- (u2);
    \draw (u3) -- (u) -- (v);
\end{tikzpicture}

\begin{tikzpicture}[scale=1]
    \node[precoloured] (u1) at (120:2cm){};
    \node[precoloured] (u2) at (90:2cm){};
    \node[precoloured] (u3) at (60:2cm){};
    \node[4vertex] (v) at (0.5,0) {};
    \node[3vertex] (u) at (-0.5,0) {};
    \node[3vertex] (w) at (1.5,0.2) {};

    \draw (u) -- (u1) -- (u2) -- (u3) -- (v) -- (u) -- (u2) -- (v);
    \draw (u3) -- (w) -- (v);
\end{tikzpicture}
\begin{tikzpicture}[scale=1]
    \node[precoloured] (u1) at (120:2cm){};
    \node[precoloured] (u2) at (90:2cm){};
    \node[precoloured] (u3) at (60:2cm){};
    \node[4vertex] (v) at (0,0) {};
    \node[3vertex] (u) at (1,0.2) {};
    \node[3vertex] (w) at (-1,0.2) {};

    \draw (v) -- (u1) -- (u2) -- (u3) -- (v) -- (u2);
    \draw (u3) -- (u) -- (v);
    \draw (u1) -- (w) -- (v);

\end{tikzpicture}}%
    \subcaptionbox{\ref{S232c}-singular}[0.4\textwidth]{
\begin{tikzpicture}[scale=1]
    \node[precoloured] (u1) at (120:2cm){};
    \node[precoloured] (u2) at (90:2cm){};
    \node[precoloured] (u3) at (60:2cm){};
    \node[4vertex] (v) at (0,0) {};
    \draw (v) -- (u1) -- (u2) -- (u3) -- (v) -- (u2);
\end{tikzpicture}


\begin{tikzpicture}[scale=1]
    \node[precoloured] (u1) at (120:2cm){};
    \node[precoloured] (u2) at (90:2cm){};
    \node[precoloured] (u3) at (60:2cm){};
    \node[3vertex] (v) at (0,0) {};
    \node[geq4vertex] (u) at (0,1) {};
    \draw (v) -- (u1) -- (u2) -- (u3) -- (v);
    \draw (u2) -- (u) -- (u1);
    \draw (v) -- (u) -- (u3);
\end{tikzpicture}

\begin{tikzpicture}[scale=1]
    \node[precoloured] (u1) at (120:2cm){};
    \node[precoloured] (u2) at (90:2cm){};
    \node[precoloured] (u3) at (60:2cm){};
    \node[4vertex] (v) at (0.5,0) {};
    \node[3vertex] (u) at (-0.5,0) {};

    \draw (u) -- (u1) -- (u2) -- (u3);
    \draw (v) -- (u) -- (u2) -- (v);
    \draw[wavy] (u3) -- (v);
\end{tikzpicture}

\begin{tikzpicture}[scale=1]
    \node[precoloured] (u1) at (120:2cm){};
    \node[precoloured] (u2) at (90:2cm){};
    \node[precoloured] (u3) at (60:2cm){};
    \node[3vertex] (v1) at (-0.6,0) {};
    \node[4vertex] (v2) at (0,0) {};
    \node[3vertex] (v3) at (0.6,0) {};
    \draw (v1) -- (u1) -- (u2) -- (u3) -- (v3) -- (v2) -- (v1);
    \draw (v1) -- (u2) -- (v2);
    \draw (u2) -- (v3);
\end{tikzpicture}

\begin{tikzpicture}[scale=1]
    \node[precoloured] (u1) at (120:2cm){};
    \node[precoloured] (u2) at (60:2cm){};
    \node[3vertex] (v) at (0,0) {};
    \draw (v) -- (u1) -- (u2) -- (v);
\end{tikzpicture}

\begin{tikzpicture}[scale=1]
    \node[precoloured] (u1) at (120:2cm){};
    \node[precoloured] (u2) at (90:2cm){};
    \node[precoloured] (u3) at (60:2cm){};
    \node[4vertex] (v) at (0,0) {};
    \node[3vertex] (u) at (1,0.2) {};
    \draw (v) -- (u1) -- (u2) -- (u3) -- (v) -- (u2);
    \draw (u3) -- (u) -- (v);
\end{tikzpicture}

\begin{tikzpicture}[scale=1]
    \node[precoloured] (u1) at (120:2cm){};
    \node[precoloured] (u2) at (90:2cm){};
    \node[precoloured] (u3) at (60:2cm){};
    \node[4vertex] (v) at (0.5,0) {};
    \node[3vertex] (u) at (-0.5,0) {};
    \node[3vertex] (w) at (1.5,0.2) {};

    \draw (u) -- (u1) -- (u2) -- (u3) -- (v) -- (u) -- (u2) -- (v);
    \draw (u3) -- (w) -- (v);
\end{tikzpicture}
\begin{tikzpicture}[scale=1]
    \node[precoloured] (u1) at (120:2cm){};
    \node[precoloured] (u2) at (90:2cm){};
    \node[precoloured] (u3) at (60:2cm){};
    \node[4vertex] (v) at (0,0) {};
    \node[3vertex] (u) at (1,0.2) {};
    \node[3vertex] (w) at (-1,0.2) {};

    \draw (v) -- (u1) -- (u2) -- (u3) -- (v) -- (u2);
    \draw (u3) -- (u) -- (v);
    \draw (u1) -- (w) -- (v);

\end{tikzpicture}}%
    \hspace*{0.02\textwidth}
\captionsetup{width=0.9\textwidth}
\caption{The singular canvases mentioned by Theorem~\ref{thm:maintech}. Vertices $v$ with $f(v)=2$ are represented by triangles, those with $f(v) = 3$ by squares.
An edge counted by $s$ is wavy.
\label{fig:singularcanvases}}
\end{figure}
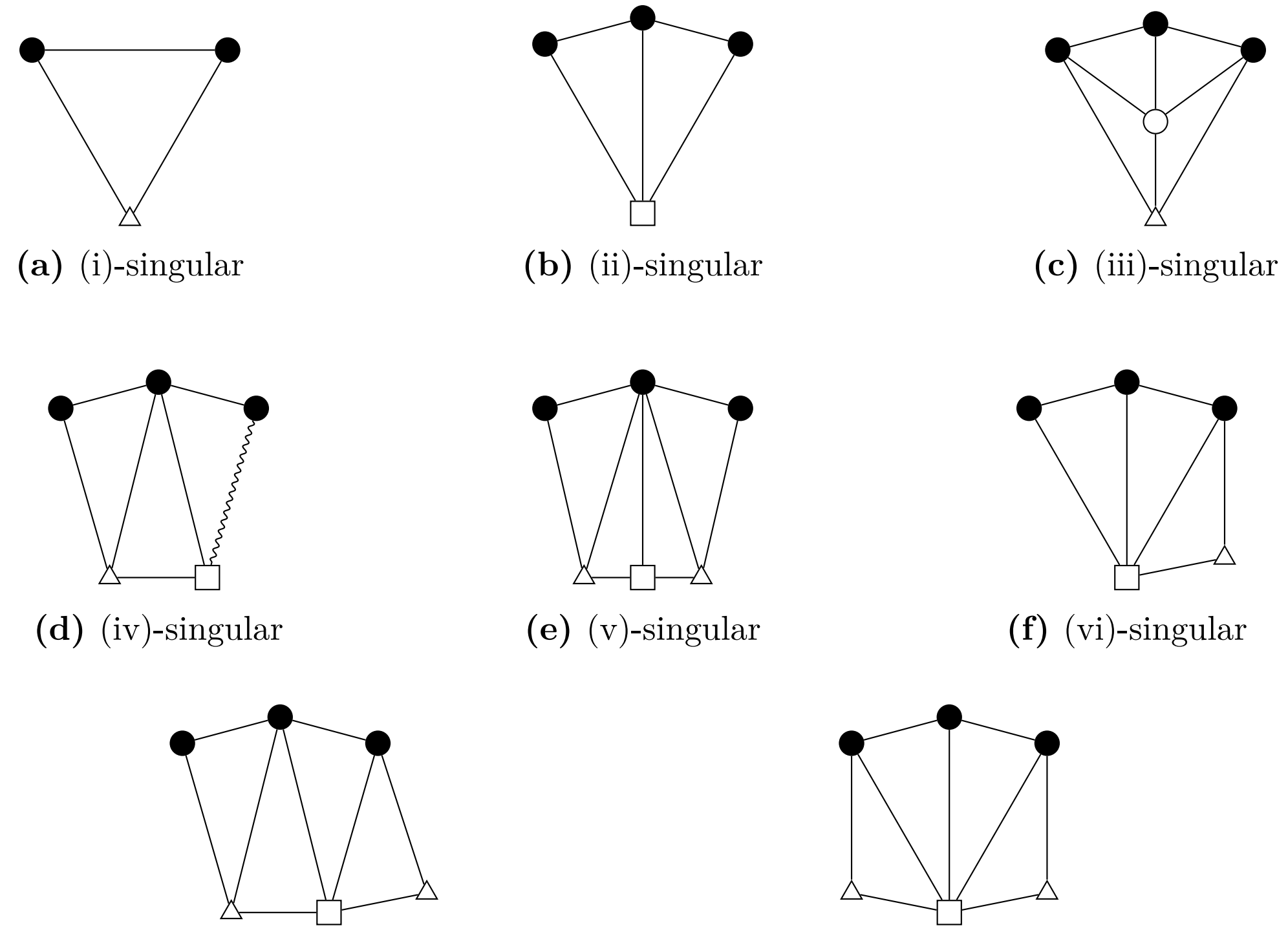

Bernshteyn, Lee, and Smith-Roberge~\cite{bernshteyn2024weak} proved that every valid tripod-free canvas is weakly degenerate.
However, excluding tripods is not sufficient to ensure $\alpha$-available weak degeneracy for any $\alpha>1$.
Indeed, the valid tripod-free canvases listed in the following definition (and illustrated in Figure~\ref{fig:singularcanvases}) are counterexamples.

\begin{definition}
    A canvas $(G,P,f)$ is \emph{singular} if $v(P)\le 2$ and $V(G)=V(P)$, or $v(P) = 2$ and
    \begin{enumerate}[(i)]
        \item\label{S} $V(G)\setminus V(P)$ consists of a single vertex $v$ with $f(v) = 2$ and $v$ is adjacent to both vertices of $P$,
    \end{enumerate}
    or $v(P)=3$ and one of the following possibilities holds:
    \begin{enumerate}[(i),resume]
        \item\label{S3} $V(G)-V(P)$ consists of a single vertex $v$ with $f(v) = 3$ and $v$ is adjacent to every vertex of $P$.
        \item\label{S2i} $G-V(P)$ is a path $uv$ with $u \in V(\partial G)$ and $v \not \in V(\partial G)$, $f(v) = 4$ and $f(u) = 2$, $v$ is adjacent to every vertex of $P$, and $u$ is adjacent to the endpoints of $P$.
        \item\label{S23a} $G-V(P)$ is a path $uv \subset \partial G$, $f(v) = 2$ and $f(u) = 3$, and each of $u$ and $v$ is adjacent to two consecutive vertices of $P$.
        \item\label{S232a} $G-V(P)$ is a path $u_1u_2u_3 \subset \partial G$, $f(u_1) = f(u_3) = 2$ and $f(u_2) = 3$, $u_1$ and $u_3$ are each adjacent to two consecutive vertices of $P$, and $u_2$ is adjacent to the middle vertex of $P$.
        \item\label{S23b} $G-V(P)$ is a path $u_1u_2 \subseteq \partial G$, $f(u_1) = 3$ and $f(u_2) = 2$, $u_1$ is adjacent to all vertices of $P$, and $u_2$ is adjacent to an endpoint of $P$.
        \item\label{S232b} $G-V(P)$ is a path $u_1u_2u_3 \subseteq \partial G$, $f(u_1) = f(u_3) = 2$ and $f(u_2) = 3$, $u_1$ and $u_2$ are each adjacent to two consecutive vertices of $P$, and $u_3$ is adjacent to the same endpoint of $P$ as $u_2$.
        \item\label{S232c} $G-V(P)$ is a path $u_1u_2u_3 \subseteq \partial G$, $f(u_1) = f(u_3) = 2$ and $f(u_2) = 3$, $u_1$ and $u_3$ are each adjacent to a different endpoint of $P$, and $u_2$ is adjacent to all three vertices of $P$.
    \end{enumerate}
 If a canvas is described by \ref{S} above, we say that it is \ref{S}-singular (similarly \ref{S3}- through \ref{S232c}-singular).
\end{definition}

For our argument, we need to measure the size of a canvas in a somewhat technical way.

\begin{definition}
    Given a canvas $K = (G,P,f)$,
    \begin{itemize}
        \item let $n(K)$ denote the number of vertices $v\in V(G)\setminus V(P)$ with $f(v) \geq 3$,
        \item let $q(K)$ denote the number of vertices $v\in V(\partial G) \setminus V(P)$ with $f(v) \geq 4$, and
        \item let $s(K)$ denote the number of edges in $\partial G$ that have exactly one endpoint in $P$ and
	      whose other endpoint $v$ satisfies
	      $$f(v)\ge \max(3, |N(v) \cap V(P)|+1).$$
              In other words, $s$ counts an edge $pv \in E(\partial G)\setminus E(P)$ with $p \in V(P)$ if
	      either $f(v) \geq 4$, or $f(v) = 3$ and $v$ has at most two neighbors in $P$. 
    \end{itemize} 
    We define
    $$\lambda(K)=n(K) + 2q(K) + 2s(K)+2\cdot(3-v(P)).$$
\end{definition}

It will be convenient to compute the values from this definition for singular canvases.
\begin{obs}\label{obs:sing}
Let $K=(G,P,f)$ be a singular canvas with $V(G)\neq V(P)$.  Then $K$ is weakly degenerate,
but every sequence of operations showing this has total availability exactly one.  Moreover,
the values of $\lambda(K)$, $n(K)$, $q(K)$, and $s(K)$ are given by the following table.
\begin{center}
\begin{tabular}{c|cccc}
$K$			& $n(K)$ & $q(K)$ & $s(K)$ & $\lambda(K)$\\
\hline
$G=P$, $v(P)=1$		& 0	 & 0	  & 0      & 4   \\
$G=P$, $v(P)=2$		& 0	 & 0	  & 0	   & 2	 \\
\ref{S}-singular	& 0      & 0      & 0      & 2   \\
\ref{S3}-singular	& 1      & 0      & 0      & 1   \\
\ref{S2i}-singular	& 1      & 0      & 0      & 1   \\
\ref{S23a}-singular	& 1      & 0      & 1      & 3   \\
\ref{S232a}-singular	& 1      & 0      & 0      & 1   \\
\ref{S23b}-singular	& 1      & 0      & 0      & 1   \\
\ref{S232b}-singular	& 1      & 0      & 0      & 1   \\
\ref{S232c}-singular	& 1      & 0      & 0      & 1   \\
\end{tabular}
\end{center}
\end{obs}

Note that the singular canvases have average availability exactly one and non-zero $\lambda$, forcing us to exclude them in the following
stronger version of our main theorem.
\begin{thm}\label{thm:maintech}
    Let $c=2^{1/\CCC}$.  For any valid tripod-free canvas $K$,
    if $K$ is non-singular, then it is totally $c^{\lambda(K)}$-availably weakly degenerate.
\end{thm}

This indeed strengthens Theorem~\ref{thm:main}.
\begin{proof}[Proof of Theorem~\ref{thm:main}, assuming Theorem~\ref{thm:maintech}]
Let $c=2^{1/\CCC}$ be the constant from the statement of Theorem~\ref{thm:maintech}.
Let $G$ be any planar graph and let $f$ be the function that assigns the value $4$ to every vertex of $G$.
We need to show that $G$ is $c$-availably weakly $f$-degenerate.  This is clearly the case if $v(G)\le 1$,
and thus suppose that $v(G)\ge 2$.

Let $p$ be any vertex incident with the outer face of $G$ and let $P$ be the single-vertex path consisting of $p$.
Then $K=(G,P,f)$ is a valid tripod-free non-singular canvas, and thus Theorem~\ref{thm:maintech} implies that $K$ is
totally $c^{\lambda(K)}$-availably weakly degenerate.  Moreover, we have $\lambda(K)\ge n(G)=v(G)-1$.  Therefore,
the graph $G-p$ is totally $c^{v(G)-1}$-availably
weakly $f_K$-degenerate.  Since the operation $\del(G,f,p)=(G-p,f_K)$ has availability $5>c$, it follows that
the graph $G$ is totally $c^{v(G)}$-availably weakly $f$-degenerate, or equivalently,
$c$-availably weakly $f$-degenerate.
\end{proof}

\section{Describing a minimum counterexample}\label{sec:properties}

Throughout the rest of the paper, let $c=2^{1/\CCC}$ refer to the constant from the statement of Theorem~\ref{thm:maintech}.
A canvas $K'=(G',P',f')$ is \emph{smaller} than another canvas $K=(G,P,f)$ if
\begin{itemize}
\item $v(G')<v(G)$, or
\item $v(G')=v(G)$ and the outer face of $G'$ has strictly \emph{larger} length than the outer face of $G$, or
\item $v(G')=v(G)$, the outer faces of $G$ and $G'$ have the same length, and $|E(G')|>|E(G)|$.
\end{itemize}
Let us emphasize that indeed the directions of inequalities in the second and third cases are as intended (even though this might seem a bit
counterintuitive).  The canvas $K$ is a \emph{counterexample} if $K$ is valid, tripod-free, non-singular, and not totally $c^{\lambda(K)}$-availably weakly degenerate.
A counterexample is \emph{minimal} if no smaller canvas is a counterexample.
Let us note the following obvious consequence of minimality, which we frequently use.
\begin{obs}\label{obs:mintolam}
Let $K=(G,P,f)$ be a minimal counterexample and let $K'=(G',P',f')$ be a canvas
obtained from $K$ by a non-empty sequence $O_1$, $O_2$, \ldots, $O_k$ of legal $\del$, $\delsave$, and $\Omit$ operations whose
availabilities are $a_1$, \ldots, $a_k$, respectively.  Let $a=\prod_{i=1}^k a_i$.  Let $K^\star=(G^\star,P^\star,f^\star)$ be a canvas such that
$v(G^\star)\le v(G')$, $G^\star-V(P^\star)=G'-V(P')$, and $f^\star_{K^\star}=f'_{K'}$.
Let $\lambda^\star=\lambda(K^\star)$ if the canvas $K^\star$ is non-singular and $\lambda^\star=0$ if it is singular.
If the canvas $K^\star$ is valid and tripod-free, then
$$\lambda(K)>\lambda^\star+\log_c a.$$
\end{obs}
\begin{proof}
Note that $v(G^\star)\le v(G')=v(G)-k<v(G)$, and thus the canvas $K^\star$ is smaller than the canvas $K$.
It follows that $K^\star$ is not a counterexample.
If the canvas $K^\star$ is non-singular, then since it is valid and tripod-free, it is totally $c^{\lambda(K^\star)}$-availably weakly degenerate.
If $K^\star$ is singular, then it is (totally $c^{0}$-availably) weakly degenerate by Observation~\ref{obs:sing}.
In either case, it follows that $K^\star$ is totally $c^{\lambda^\star}$-availably weakly degenerate.

Since $G^\star-V(P^\star)=G'-V(P')$ and $f^\star_{K^\star}=f'_{K'}$, the canvas $K'$ is also totally $c^{\lambda^\star}$-availably weakly degenerate.  Therefore, the canvas $K$ is
totally $(a\cdot c^{\lambda^\star})$-availably weakly degenerate.  Since $K$ is a counterexample, it is
not totally $c^{\lambda(K)}$-availably weakly degenerate, and thus
$$c^{\lambda(K)}>a\cdot c^{\lambda^\star}.$$
This is equivalent to
$$\lambda(K)>\lambda^\star+\log_c a.$$
\end{proof}

Let us remark that we usually apply this observation with $K^\star=K'$; that is, we perform a sequence of operations
on a minimal counterexample so that the resulting canvas $K'$ is valid and tripod-free, and then use the inequality from
the statement of Observation~\ref{obs:mintolam} to obtain a contradiction.   However, occasionally the canvas $K'$ itself
is not valid, but we can make it valid by adding an edge between a vertex in $P'$ and a vertex $v$ outside of $P'$
and increasing the value of $f'(v)$, which affects neither the graph $G'-V(P')$ nor the function $f_{K'}$.
In this case, the resulting canvas plays the role of $K^\star$ from Observation~\ref{obs:mintolam}.

In the rest of this section, we establishing basic properties of each minimal counterexample $K=(G,P,f)$:
Namely, we show that the graph $G$ is 2-connected
(Lemma~\ref{lemma:2conn}), that its triangles bound internal faces (Lemma~\ref{lemma:noseptriang}),
and we restrict its short generalized chords of length at most two (Lemmas~\ref{lemma:1chords} and \ref{lemma:2chords}).
In addition, we show that the $f$-value of certain vertices in $V(\partial G)$ is not too low compared to their degree
(Lemmas~\ref{lemma:minfvals-A} and \ref{lemma:minfvalues-notA}).  Let us note the following observation which
we apply in several of these lemmas.
\begin{obs}\label{obs:comb}
Let $K=(G,P,f)$ be a minimal counterexample and let $Q$ be a generalized chord of length at most two separating $G$ into subgraphs $G_1$ and $G_2$.
Suppose that $P\subseteq G_1$, and let $K_1=(G_1,P,f)$ and $K_2=(G_2,Q,f)$.
If $Q$ is an induced path in $G_2$ and the canvas $K_2$ is tripod-free, then either
\begin{enumerate}[(a)]
\item\label{comb:singS} $K_1$ is singular, $v(Q)=2$ and $K_2$ is \ref{S}-singular, and $\lambda(K)>0$, or
\item\label{comb:nonS23a} $K_1$ is non-singular, $v(Q)=3$ and $K_2$ is \ref{S23a}-singular, and $s(K)=s(K_1)+1$, or
\item\label{comb:S23anon} $v(P)=3$ and $K_1$ is \ref{S23a}-singular, $v(Q)=2$ and $K_2$ is non-singular, and $s(K)=s(K_2)+1$.
\end{enumerate}
\end{obs}
\begin{proof}
Note that the definition of a generalized chord implies that $G_1\neq P$ and $G_2\neq Q$.
For $i\in\{1,2\}$, let $\lambda_i=\lambda(K_i)$ if the canvas $K_i$ is non-singular and $\lambda_i=0$ otherwise.
Since $Q$ is an induced path in $G_2$ and since $K$ is a valid canvas, both $K_1$ and $K_2$ are valid canvases.
Since $K$ is tripod-free, $K_1$ is tripod-free as well (and $K_2$ is tripod-free by assumption).
For $i\in \{1,2\}$, the canvas $K_i$ is not a counterexample by the minimality of $K$, and thus if $K_i$ is non-singular,
then it is totally $c^{\lambda_i}$-availably weakly degenerate.  This is also the case if $K_i$ is singular, since
then $\lambda_i=0$ and the canvas $K_i$ is weakly degenerate by Observation~\ref{obs:sing}.  Let $S_i$ be an $\Omit$-free sequence of operations
showing that $K_i$ is totally $c^{\lambda_i}$-availably weakly degenerate.  

By concatenating the sequences $S_1$ and $S_2$, we see that the canvas $K$ is totally $c^{\lambda_1+\lambda_2}$-availably weakly degenerate.
Since $K$ is a counterexample and thus not totally $c^{\lambda(K)}$-availably weakly degenerate, it follows that
$$\lambda(K)>\lambda_1+\lambda_2.$$
Let $q=0$ if $v(Q)\le 2$ and $q=1$ if $v(Q)=3$; let us note that $q+3-v(Q)\ge 1$. Observe that $n(K_1)+n(K_2)=n(K)$ and $q(K_1)+q(K_2)=q(K)+q$,
since if $v(Q)=3$, then the middle vertex of $Q$ contributes to $q(K_1)$, but not to $q(K)$.

Furthermore, note that $s(K_1)+s(K_2)\ge s(K)$: Consider any edge $pv$ that contributes to $s(K)$.
If $v\in V(G_1)$, then this edge also contributes to $s(K_1)$.  If $v\in V(G_2)\setminus V(G_1)$,
then $p\in V(Q)$, and $pv$ contributes to $s(K_2)$ unless $f(v)=3$ and $v$ has three neighbors in $Q$.
However, in this case we have $v(Q)=3$ and we get an additional contribution to $s(K_1)$ for the edge $pu$, where $u$ is the middle vertex of $Q$.

Therefore,
\begin{align}
\lambda(K_1)+\lambda(K_2)&=(n(K_1)+2q(K_2)+2s(K_1)+2\cdot(3-v(P))\nonumber\\
&\phantom{{}={}}+(n(K_2)+2q(K_2)+2s(K_2)+2\cdot(3-v(Q)))\nonumber\\
&=n(K)+2q(K)+2q+2s(K)+2\cdot(s(K_1)+s(K_2)-s(K))\nonumber\\
&\phantom{{}={}}+2\cdot(3-v(P))+2\cdot(3-v(Q))\nonumber\\
&=\lambda(K)+2q+2\cdot(3-v(Q))+2\cdot(s(K_1)+s(K_2)-s(K))\nonumber\\
&>\lambda_1+\lambda_2+2\cdot(q+3-v(Q))+2\cdot (s(K_1)+s(K_2)-s(K))\label{eq:comb}\\
&\ge \lambda_1+\lambda_2+2.\nonumber
\end{align}
It follows that $\lambda(K_1)\neq \lambda_1$ or $\lambda(K_2)\neq \lambda_2$ and thus $K_1$ or $K_2$ is singular.

If both $K_1$ and $K_2$ are singular, then $\lambda(K)>\lambda_1+\lambda_2=0$.  Moreover,
we cannot have $v(Q)=3$, since in that case we would have $q(K_1)\ge q>0$, contradicting Observation~\ref{obs:sing}.
It follows that $K_2$ is $\ref{S}$-singular and we obtain the first outcome of the lemma.

Therefore, we may also assume that there exists $i\in \{1,2\}$ such that $K_i$ is non-singular.
Then $\lambda(K_i)=\lambda_i$, $K_{3-i}$ is singular, and $\lambda_{3-i}=0$.
Thus, (\ref{eq:comb}) gives
$$\lambda(K_{3-i})>2\cdot(q+3-v(Q))+2\cdot(s(K_1)+s(K_2)-s(K))\ge 2+2\cdot (s(K_1)+s(K_2)-s(K))\ge 2.$$
By Observation~\ref{obs:sing}, this is only possible if $K_{3-i}$ is \ref{S23a}-singular
and $s(K)=s(K_1)+s(K_2)$. Moreover, since $K_{3-i}$ is \ref{S23a}-singular, we have $s(K_{3-i})=1$,
and thus $s(K)=s(K_i)+1$.  Note also that $q+3-v(Q)=1$, and thus $v(Q)\in\{2,3\}$;
and if $v(Q)=3$, then $1=q\le q(K_1)$.  Since $K_{3-i}$ is \ref{S23a}-singular, it satisfies $q(K_{3-i})=0$,
and thus if $i=2$, then $v(Q)=2$.
Thus, we obtain the second outcome of the lemma when $i=1$ and the third one when $i=2$.
\end{proof}

Moreover, let us note the following non-triviality observation.
\begin{obs}\label{obs:ge2}
    If $K=(G,P,f)$ is a counterexample, then $v(G)\ge v(P)+2$.
\end{obs}
\begin{proof}
Suppose first that $v(G)=v(P)$.  Since $K$ is valid, $P$ is an induced path in $G$, and thus $G=P$.
Since $K$ is non-singular, we have $v(P)=3$.  However, then $\lambda(K)=0$, and since $K$ is (totally $c^0$-availably)
weakly degenerate, this contradicts the assumption that $K$ is a counterexample.

Suppose now for a contradiction that $v(G)=v(P)+1$, and thus $V(G)\setminus V(P)$ consists of a single vertex $v$.
We have $f_K(v)\ge 1$: This follows from the validity of $K$ if $v(P)=1$, from validity and non-\ref{S}-singularity if $v(P)=2$,
and from validity, tripod-freeness, and non-\ref{S3}-singularity if $v(P)=3$.
Thus, the operation $\del(K,v)$ has availability at least two,
and thus $K$ is totally $2$-availably weakly degenerate.  Since $K$ is a counterexample and $2=c^7$,
it follows that $\lambda(K)>\CCC$.  Note that $n(K)\le 1$, $q(K)\le 1$, and $s(K)\le 2$,
and thus $\lambda(K)\le 7+2\cdot(3-v(P))$, which implies that $v(P)\le 2$.

Moreover, note that $3-v(P)+s(K)\le 3$, since either $v(P)=2$ and $s(K)\le 2$, or $v(P)=1$ and $s(K)\le 1$.
Hence, $\lambda(K)\le 7+2q(K)$, and thus $q(K)>0$.  Hence, the vertex $v$ contributes to $q(K)$, and thus $f(v)\ge 4$.
Since $v(P)\le 2$, this implies that $f_K(v)\ge 2$, and thus $K$ is actually totally $3$-availably weakly degenerate.
However, since $\lambda(K)\le 7+2q(K)\le 9$ and $3>c^9$, this contradicts the assumption
that $K$ is a counterexample.
\end{proof}

It is now easy to establish a lower bound on the connectivity of minimal counterexamples.

\begin{lemma}\label{lemma:2conn}
    If $K=(G,P,f)$ is a minimal counterexample, then the graph $G$ is connected and has no $0$-chords.
\end{lemma}
\begin{proof}
Suppose first that $G$ is not connected.  Let $M$ be the component of $G$ containing $P$, and note that there exists a component $R\neq M$ of $G$
such that a face $g$ of $G$ is incident with both a vertex $u\in V(M)$ and a vertex $v\in V(R)$.
Observe furthermore that we can choose $u$ so that either $u\in V(P)$, or
$u\not\in V(P)$ and $f(u)\ge 3$, since $f^{-1}(2)\cap V(G)\setminus V(P)$ is an independent set. 
Consequently, the canvas $K'=(G+uv,P,f)$ is valid.  Moreover, since $uv$ is a bridge in $K'$ and since $K$ is tripod-free,
it follows that $K'$ is tripod-free and non-singular.  Since $K'$ has $v(K)$ vertices,
at least as long outer face as $K$, and more edges than $K$, it is smaller than $K$, and thus it is not a counterexample.
Thus, $K'$ is totally $c^{\lambda(K')}$-availably weakly degenerate,
and Observation~\ref{obs:monot} implies that $K$ is totally $c^{\lambda(K')}$-availably weakly degenerate as well.
Moreover, we have $\lambda(K')\ge \lambda(K)$ (adding the edge $uv$ can only increase the $s$-value),
and thus $K$ is also totally $c^{\lambda(K)}$-availably weakly degenerate.
This is a contradiction.

Therefore, $G$ is connected.  Suppose now that there exists a 0-chord $Q$ of $G$ that separates $G$
into subgraphs $G_1$ and $G_2$, with the labels of $G_1$ and $G_2$ chosen so that $v(G_1 \cap P) \geq v(G_2 \cap P)$.
Let $K_1 = (G_1, P \cap G_1, f)$ and let $v$ denote the single vertex of $Q$.

If $P\subseteq G_1$, then let $K_2 = (G_2, Q, f)$.
Note since $v(Q)=1$, the canvas $K_2$ is non-singular and tripod-free.
However, then we obtain contradiction by Observation~\ref{obs:comb}.

Therefore, $P \not \subseteq G_1$.  Note that this is only possible if $P$ is a path of length two,
$v$ is the middle vertex of $P$, and $v(G_1 \cap P)=v(G_2 \cap P)=2$.  Let $K_2 = (G_2, P \cap G_2, f)$.
Clearly both canvases $K_1$ and $K_2$ are valid and tripod-free.
For $i\in \{1,2\}$, if $K_i$ is non-singular, then let $\lambda_i=\lambda(K_i)$, otherwise let $\lambda_i=0$.
Hence, the minimality of the counterexample $K$ and Observation~\ref{obs:sing} implies that the canvas $K_i$ is
totally $c^{\lambda_i}$-availably weakly degenerate.  By concatenating the two sequences of operations showing
this, we see that $K$ is totally $c^{\lambda_1+\lambda_2}$-availably weakly degenerate.
Since $K$ is a counterexample, it follows that
$$\lambda(K)>\lambda_1+\lambda_2.$$
On the other hand, we have $n(K)=n_1(K)+n_2(K)$, $q(K)=q(K_1)+q(K_2)$, and $s(K)=s(K_1)+s(K_2)$,
and thus
\begin{align*}
\lambda(K)&=\lambda(K)-2\cdot(3-v(P))=n(K) + 2q(K) + 2s(K)\\
&=(n(K_1) + 2q(K_1) + 2s(K_1))+(n(K_2) + 2q(K_2) + 2s(K_2))\\
&=(\lambda(K_1)-2\cdot(3-v(P\cap G_1)))+(\lambda(K_2)-2\cdot(3-v(P\cap G_2)))\\
&=\lambda(K_1)+\lambda(K_2)-4.
\end{align*}
Therefore,
$$\lambda(K_1)+\lambda(K_2)>\lambda_1+\lambda_2+4.$$
By symmetry, we can assume that $\lambda(K_1)>\lambda_1+2$.
However, by the definition of $\lambda_1$, we either have $\lambda(K_1)=\lambda_1$,
or $K_1$ is singular and Observation~\ref{obs:sing} (taking into account that $v(P\cap G_1)=2$) gives $\lambda(K_1)=2$.
This is a contradiction.
\end{proof}

In particular, Lemma~\ref{lemma:2conn} implies that $\partial G$ is a cycle.
A triangle $T$ in a plane graph $G$ is \emph{empty} if the internal face $g$ of the graph $T$
(whose drawing is inherited from $G$) is also a face of $G$, and \emph{non-empty} otherwise (in this case, at least
one vertex of $G$ is drawn in $g$).

\begin{lemma}\label{lemma:noseptriang}
    If $K=(G,P,f)$ is a minimal counterexample, then all triangles in the graph $G$ are empty.
\end{lemma}
\begin{proof}
    Suppose for a contradiction that $T=u_1u_2u_3$ is a non-empty triangle in $G$.
    Let $X$ be the set of vertices of $G$ drawn in the internal face of the graph $T$ (whose drawing is inherited from $G$).
    Let $G_1=G-X$, $G_2=G[V(T)\cup X]-u_1u_3$, and $P_2=u_1u_2u_3$.
    Let us consider the canvases $K_1=(G_1,P,f)$ and $K_2=(G_2,P_2,f)$.
    Note that both are valid.  Moreover, $K_1$ is tripod-free since $K$ is tripod-free, and $K_2$ is tripod-free and non-singular
    since $f(v)\ge 4$ holds for every vertex $v\in V(G_2)\setminus V(P_2)=X$.

    Clearly $v(K_1)<v(K)$, since $X\neq\emptyset$.  We could have $v(K_2)=v(K)$, but in that case $T$ bounds the
    outer face of $G$ and the graph $G_2$ has strictly longer outer face.  Hence, the canvases $K_1$ and $K_2$
    are smaller than $K$, and thus they are not counterexamples.

    Let $\lambda_1=\lambda(K_1)$ if $K_1$ is non-singular and $\lambda_1=0$ otherwise, so that
    $K_1$ is totally $c^{\lambda_1}$-availably weakly degenerate.  Clearly $K_2$ is $c^{\lambda(K_2)}$-availably weakly degenerate.
    By concatenating the $\Omit$-free sequences witnessing this, we see that $K$ is totally $c^{\lambda_1+\lambda(K_2)}$-availably weakly degenerate.
    Since $K$ is a counterexample, we get
    $$\lambda(K)>\lambda_1+\lambda(K_2).$$
    Note that $n(K)=n(K_1)+n(K_2)$, $q(K)=q(K_1)$ and $q(K_2)\ge 1$ (since at least one vertex of $X$ is incident with the outer face of $G_2$),
    and $s(K)=s(K_1)$ and $s(K_2)\ge 1$ (the connectivity of $G$ implies that there is at least one edge between $X$ and $V(P_2)$
    contained in the boundary of the outer face of $G_2$).  Consequently,
    \begin{align*}
    \lambda(K)&=n(K)+2q(K)+2s(K)+2\cdot(3-v(P))\\
              &=n(K_1)+n(K_2)+2q(K_1)+2s(K_1)+2\cdot(3-v(P))=\lambda(K_1)+n(K_2)\\
	      &=\lambda(K_1)+\lambda(K_2)-2q(K_2)-2s(K_2)-2\cdot (3-v(P_2))\\
	      &\le \lambda(K_1)+\lambda(K_2)-4.
    \end{align*}
    Therefore, we get
    $$\lambda(K_1)\ge \lambda(K)-\lambda(K_2)+4>\lambda_1+4.$$
    However, this contradicts the definition of $\lambda_1$ and Observation~\ref{obs:sing}.
\end{proof} 

The following lemma will greatly simplify further arguments.
\begin{lemma}\label{lemma:minfvals-A}
    Let $K=(G,P,f)$ be a minimal counterexample and let $v\in V(G)\setminus V(P)$ be a vertex incident with the outer
    face of $G$.  If $f(v)=2$, then $\deg v\geq 3$.
\end{lemma}
\begin{proof}
Suppose for a contradiction that $\deg v\le 2= f(v)$, and thus we can legally $\Omit$ the vertex $v$ with availability
$3-\deg v\ge 1$.
Consider the canvas $K'=(G-v,P,f)=\Omit(K,v)$.  Since $K$ is valid and tripod-free, the canvas $K'$ is valid and tripod-free as well.
Note that $$\lambda(K')\ge\lambda(K),$$ since the vertices with $f$-value two do not contribute to $\lambda$.
Let $\lambda'=\lambda(K')$ if $K'$ is non-singular and $\lambda'=0$ otherwise.  Observation~\ref{obs:mintolam} (with $K^\star=K'$)
implies that 
$$\lambda(K)>\lambda' + \log_c (3-\deg v)\ge \lambda'.$$
Therefore, we have $\lambda(K')>\lambda'$, and thus $K'$ is singular, $\lambda'=0$, and
$$\lambda(K')>\log_c(3-\deg v).$$
Observation~\ref{obs:sing} gives $\lambda(K')\le 4$, and since $\log_c 2=\CCC>4$, it follows
that $\deg v=2$.  By Observation~\ref{obs:ge2}, we have $G-v\neq P$.
Since the canvas $K$ is valid, an inspection of the singular canvases shows that two neighbors of $v$
are adjacent in $\partial (G-v)-E(P)$, the canvas $K'$ is \ref{S3}-, \ref{S23a}-, or \ref{S23b}-singular,
and one of the neighbors of $v$ is an end of $P$.
However, then $K$ is \ref{S23b}-, \ref{S232b}-, or \ref{S232c}-singular, respectively,
contradicting the assumption that $K$ is a counterexample.
\end{proof}

Recall that by Lemma \ref{lemma:2conn}, $\partial G$ is a cycle. Next, we show that this cycle is induced.

\begin{lemma}\label{lemma:1chords}
    If $K=(G,P,f)$ is a minimal counterexample, then $G$ does not have a $1$-chord.
\end{lemma}
\begin{proof}
    Suppose for a contradiction that $Q=u_1u_2$ is a 1-chord separating $G$ into subgraphs $G_1$ and $G_2$, where
    $v(P \cap G_1) \geq v(P \cap G_2)$.  Choose such a 1-chord $Q$ with $v(P\cap G_1)$ maximal.
    Let us now distinguish two cases.

    \paragraph{Case 1: $P\subseteq G_1$.} Consider the canvases $K_1=(G_1,P,f)$ and $K_2=(G_2,Q,f)$.  
    Since $v(Q)=2$, the canvas $K_2$ is tripod-free,
    and since $Q$ is a $1$-chord, we have $G_1\neq P$ and $G_2\neq Q$.  Lemma~\ref{lemma:minfvals-A}
    implies that $K_2$ is not \ref{S}-singular. Therefore, again since $v(Q) = 2$,  Observation~\ref{obs:comb} implies that  
    $v(P)=3$ and $K_1$ is \ref{S23a}-singular,  and $K_2$ is non-singular, and $s(K)=s(K_2)+1$.
    Let $\partial G_1=p_3p_2p_1v_1v_2$, where $P=p_1p_2p_3$, $f(v_1)=2$ and $f(v_2)=3$.

    Note that $Q\neq v_2p_3$, since in that case the edge $v_2p_3$ would not be contributing to $s(K)$
    and we would have $s(K)\le s(K_2)$.  Hence, we have $Q\in\{p_1v_1,v_1v_2\}$.  Let $u$ be the neighbor
    of $v_1$ in $\partial G_2-E(Q)$.  Since $K$ is valid and $f(v_1)=2$, we have $f(u)\ge 3$.
    Therefore, the edge $v_1u$ contributes to $s(K_2)$, but not $s(K)$.  Since the only edge of $K_1$
    contributing to $s(K)$ is $v_2p_3$, it follows that
    $s(K)\le s(K_2)$.  This is a contradiction, since $s(K)=s(K_2)+1$.

    \paragraph{Case 2: $P\not\subseteq G_1$.}  Since $v(P \cap G_1) \geq v(P \cap G_2)$,
    it follows that $v(P)=3$ and one endpoint of $Q$ (say $u_2$) is the middle vertex of $P$.
    We can label the vertices $P$ as $p_1pp_2$ so that $p_1\in V(G_1)$, $p=u_2$, and $p_2\in V(G_2)$.
    By the maximality from the choice of $Q$, we see that
    \begin{enumerate}[($\dag$)]
    \item\label{prop:allwithmiddle} all $1$-chords of $G$ are incident with $p$.
    \end{enumerate}
    Moreover, we can without loss of generality assume that if $p$ has a neighbor not in $P$ with $f$-value $2$,
    then $u_1$ was chosen as one of them; that is,
    \begin{enumerate}[($\star$)]
    \item\label{prop:twoifcan} either $f(u_1)=2$, or every neighbor $u\in V(G)\setminus V(P)$ of $p$ satisfies
    $f(u)\ge 3$.
    \end{enumerate}
    If $u_1$ were adjacent both to $p_1$ and $p_2$, then \ref{prop:allwithmiddle} and Lemma~\ref{lemma:noseptriang}
    would imply that $V(G)\setminus V(P)=\{u_1\}$, contradicting Observation~\ref{obs:ge2}.
    Therefore, we can by symmetry assume that $u_1$ is not adjacent to $p_2$.

    Let $P_1=P\cap G_1$ and $K_1=(G_1,P_1,f)$, let $P_2$ be the path of length two in $G_2$ formed by $P\cap G_2$ and $Q$, and
    let $K_2=(G_2,P_2,f)$.  Since the canvas $K$ is valid, both $K_1$ and $K_2$ are valid as well (the path $P_1$ has length one and the
    path $P_2$ is induced, since $u_1$ is not adjacent to $p_2$).  Moreover, both are tripod-free:
    For $K_1$ this is trivial since $v(P_1)=2$, and if $K_2$ contained a tripod subcanvas,
    then $G_2$ would contain a triangle $pu_1w$ with $f(w)=2$, and \ref{prop:twoifcan} would imply that $f(u_1)=2$,
    contradicting the validity of $K$.

    For $i\in\{1,2\}$, let $\lambda_i=\lambda(K_i)$ if $K_i$ is non-singular and $\lambda_i=0$ otherwise,
    so that minimality of $K$ and Observation~\ref{obs:sing} imply that $K_i$ is totally $c^{\lambda_i}$-availably weakly degenerate.
    By concatenating the $\Omit$-free sequences of operations that witness this, we see
    that $K$ is totally $c^{\lambda_1+\lambda_2}$-availably weakly degenerate (note that here it is important that $u_1$ is not
    adjacent to $p_2$, since this ensures that $f_K(u_1)=f_{K_1}(u_1)$).
    Since $K$ is a counterexample, we get
    $$\lambda(K)>\lambda_1+\lambda_2.$$
    Observe that $n(K)=n(K_1)+n(K_2)$ and $q(K)=q(K_1)+q(K_2)$.  Each edge of $G$ that contributes to $s(K)$ also contributes
    to $s(K_1)$ or to $s(K_2)$.  Moreover,
    \begin{itemize}
    \item in the case that $f(u_1)\ge 3$, the edge $pu_1$ contributes to $s(K_1)$; and,
    \item in the case that $f(u_1)=2$, consider the neighbor $w$ of $u_1$ in $\partial G_2-p$; by the validity of $K$, we have $f(w)\ge 3$,
    and thus the edge $u_1w$ contributes to $s(K_2)$, unless $w$ is adjacent to all vertices of $P_2$ and $f(w)=3$.
    \end{itemize}
    Hence, we have $s(K_1)+s(K_2)\ge s(K)+s$, where $s=0$ if $f(u_1)=2$, $f(w)=3$, and $w$ is adjacent to $p$ and $p_2$,
    and $s=1$ otherwise.  Therefore,
    \begin{align*}
    \lambda(K)&=n(K)+2q(K)+2s(K)+2\cdot(3-v(P))\\
              &\le n(K_1)+n(K_2)+2q(K_1)+2q(K_2)+2s(K_1)+2s(K_2)-2s+2\cdot(3-v(P_2))\\
	      &=\lambda(K_1)+\lambda(K_2)-2s-2\cdot(3-v(P_1))\\
	      &=\lambda(K_1)+\lambda(K_2)-2s-2.
    \end{align*}
    It follows that
    \begin{equation}\label{eq:1chord}
    \lambda(K_1)+\lambda(K_2)>\lambda_1+\lambda_2+2s+2.
    \end{equation}
    This implies that there exists $i\in\{1,2\}$ such that $\lambda(K_i)>\lambda_i+s+1\ge 1$.
    Since $\lambda(K_i)>\lambda_i$, the canvas $K_i$ is singular.
    Note that $G_i\neq P_i$, and since $\lambda(K_i)>1$, Observation~\ref{obs:sing}
    implies that $K_i$ is \ref{S}- or \ref{S23a}-singular.

    If $K_i$ is \ref{S}-singular, then $v(P_i)=2$, and thus $i=1$ and $f(u_1)=2$.
    Since $2=\lambda(K_1)>s+1$, we have $s=0$, and thus the neighbor $w$ of $u_1$ in $\partial G_2-p$
    satisfies $f(w)=3$ and $w$ is adjacent to $p$ and $p_2$.  But then \ref{prop:allwithmiddle} and Lemma~\ref{lemma:noseptriang}
    imply that $V(G)\setminus V(P)=\{u_1,w\}$ and $K$ is \ref{S23a}-singular, which is a contradiction.

    Therefore, $K_i$ is \ref{S23a}-singular; and in particular $v(P_i)=3$, and thus $i=2$.
    Observe that in this case $s=1$, since no vertex of $K_2$ is adjacent to all vertices of $P_2$, and that $\lambda(K_2)=3$.
    Therefore, (\ref{eq:1chord}) gives
    $$\lambda(K_1)>\lambda_1+1.$$
    It follows that $K_1$ is singular, and since $v(P_1)=2$ and $G_1\neq P_1$, the canvas $K_1$ is actually \ref{S}-singular.  However, then
    $K$ is \ref{S232b}-singular, which is a contradiction.
\end{proof}

A plane graph $G$ is a \emph{near-triangulation} if all internal faces of $G$ have length three.
\begin{lemma}\label{lemma:neartri}
If $K=(G,P,f)$ is a minimal counterexample, then $G$ is a near-triangulation, $G$ has minimum degree at least three,
and all vertices of $G$ not incident with the outer face have degree at least four.
\end{lemma}
\begin{proof}
Suppose for a contradiction that an internal face $g$ of $G$ has length greater than three.
Note that there exists a cycle $C$ in $G$ such that
\begin{itemize}
\item all vertices and edges of $C$ are all incident with $g$, and
\item the face $g$ is a subset of the internal face $g'$ of the graph $C$ (whose drawing is inherited from $G$).
\end{itemize}
Let us remark that $g'$ is not necessarily equal to $g$, since at this moment we do not know (yet)
that the graph $G$ is $2$-connected.  Lemma~\ref{lemma:noseptriang} implies that the cycle $C$ has length at least
four ($C$ cannot be an empty triangle because of the face $g$ of length at least four
contained in $g'$).  Let us now distinguish two cases.

\paragraph{Case 1: There exists a vertex $v$ incident with $g$ that is not contained in $V(\partial G)$.}
Let us remark that the vertex $v$ may or may not be contained in $C$.
Let $H$ be the subgraph of $G$ consisting of $C$ and all edges of $G$ between $v$ and vertices in $V(C)\setminus\{v\}$.
\begin{itemize}
\item If $v\not\in V(C)$, then the edges of $E(H)\setminus E(C)$ are drawn in $g'$, and since all edges of $C$ are incident with $g$,
there is at most one such edge.
\item If $v\in V(C)$, then the edges of $E(H)\setminus E(C)$ are drawn outside of $g'$, since all edges of $C$ are incident with $g$.
In this case Lemma~\ref{lemma:noseptriang} implies that the outer face of $H$ (with the drawing inherited from $G$) has length at least four.
\end{itemize}
In either case, observe that this implies that there exists a vertex $u\in V(C)\setminus \{v\}$ non-adjacent to $v$.
Consider the canvas $K'=(G+uv,P,f)$, with the edge $uv$ drawn inside $g$.
Since $f(v)\ge 4$ and $K$ is valid and tripod-free, the canvas $K'$ is also valid and tripod-free.
Since $G+uv$ has the same number of vertices and the same outer face as $G$, but more edges,
the canvas $K'$ is smaller than $K$, and thus $K'$ is not a counterexample.

If $K'$ is non-singular, then it is $c^{\lambda(K')}$-availably weakly degenerate, and 
Observation~\ref{obs:monot} implies that $K$ is totally $c^{\lambda(K')}$-availably weakly degenerate as well.
However, note that $\lambda(K)=\lambda(K')$, and thus this contradicts the assumption that $K$ is a counterexample.

Therefore, the canvas $K'$ is singular.  Since $f(v)\ge 4$, this is only possible if $f(v)=4$ and $K'$ is \ref{S2i}-singular.
Let $z$ denote the vertex in $V(G)\setminus (V(P)\cup \{v\})$, where $f(z)=2$.
Note that $\lambda(K)=\lambda(K')=1$, and that $f_{K'}(v)=1$ and either
\begin{itemize}
\item $u=z$ and $f_K(v)=1$, or
\item $u\neq z$ and $f_K(v)=2$.
\end{itemize}
Therefore, the sequence of legal operations $\del(\del(K,z),v)$ has total availability $2>c^{\lambda(K)}$,
contradicting the assumption that $K$ is a counterexample.

\paragraph{Case 2: All vertices incident with $g$ are contained in $V(\partial G)$.}
In particular, we have $V(C)\subseteq V(\partial G)$.  Lemma~\ref{lemma:1chords} then implies that $C=\partial G$ and that $C$ is an induced
cycle.  It also follows that all vertices incident with $g$ are contained in $V(C)$, and thus the cycle
$C$ bounds the internal face $g$ of $G$.  Therefore, we have $G=C$.  Let $v$ be a vertex in $V(C)\setminus V(P)$ and
let $p\in V(P)$ be a vertex non-adjacent to $v$.  Since $\deg v=2$, Lemma~\ref{lemma:minfvals-A} implies that $f(v)\ge 3$.
Let us again define $K'=(G+uv,P,f)$, and note that since $f(v)\ge 3$ and the canvas $K$ is valid and tripod-free, the canvas $K'$ is also
valid and tripod-free.  As in the previous case, we conclude that since $K$ is a minimal counterexample and $\lambda(K')=\lambda(K)$,
the canvas $K'$ has to be singular.  Since $G$ is a cycle, $K'$ is \ref{S3}-singular.
However, then $v(G)=v(P)+1$, contradicting Observation~\ref{obs:ge2}.

\bigskip

Since both cases led to a contradiction, we conclude that $G$ is a
near-triangulation.  Consider any vertex $y\in V(G)$.  We need to argue that $\deg y\ge 3$, and
that if $y$ is not incident with the outer face, then $\deg y\ge 4$.

If $y$ is not incident with the outer face, then since $G$ is a near-triangulation, the neighbors of $y$ form a cycle
of length $\deg y$. Lemma~\ref{lemma:noseptriang} then implies the desired inequality $\deg y \geq 4$.

Finally, suppose for a contradiction that $y$ is incident with the outer face and $\deg y=2$.
Hence, the two neighbors $y_1$ and $y_2$ of $y$ both belong to $\partial G$.  Since $G$ is a near-triangulation,
the vertices $y_1$ and $y_2$ are adjacent.
Lemma~\ref{lemma:1chords} implies that $y_1y_2$ is not a $1$-chord of $G$,
and thus $\partial G$ is the triangle $yy_1y_2$.  Lemma~\ref{lemma:noseptriang} implies that actually $G$ is equal to this triangle,
and Observation~\ref{obs:ge2} gives $v(P)=1$.  Without loss of generality, we can assume that $V(P)=\{y\}$,
and Lemma~\ref{lemma:minfvals-A} implies that $f(y_1),f(y_2)\ge 3$.  Therefore, we have $f_K(y_1),f_K(y_2)\ge 2$.  However,
then the sequence of legal operations $\del(\del(K,y_1),y_2)$ has total availability at least $6$.
Observe that $n(K)=2$, $q(K)\le 2$, and $s(K)=2$, and thus $\lambda(K)\le 14$.
However then $c^{\lambda(K)}=4<6$, and thus $K$ is not a counterexample, which is a contradiction.
\end{proof}

We now prove a lemma similar to Lemma \ref{lemma:minfvals-A}, which will shorten the analysis in subsequent arguments.

\begin{lemma}\label{lemma:minfvalues-notA}
Let $K=(G,P,f)$ be a minimal counterexample and let $v\in V(G)\setminus V(P)$ be a vertex incident with the outer
face of $G$.  If $f(v)=3$, then $\deg v\geq 4$.
\end{lemma}
\begin{proof}
Let $v_1$ and $v_2$ be the two neighbors of $v$ in the cycle $\partial G$.
Suppose for a contradiction that $\deg v\le 3$; Lemma~\ref{lemma:neartri} implies that $\deg v=3$.

Note that $\deg v=f(v)$, and thus we can legally $\Omit$ the vertex $v$ (with availability one).
Consider the canvas $K'=(G-v,P,f)=\Omit(K,v)$.  Since $K$ is valid and tripod-free, the canvas $K'$ is valid and tripod-free as well.
Observation~\ref{obs:ge2} implies that $G-v\neq P$.

Let $u$ denote the neighbor of $v$ different from $v_1$ and $v_2$.
Lemma~\ref{lemma:1chords} implies that $u$ is not incident with the outer face of $G$, and thus $f(u)\ge 4$.
Since $u$ is incident with the outer face of $G-v$, the canvas $K'$ is not singular.
Hence, Observation~\ref{obs:mintolam} gives
$$\lambda(K)>\lambda(K').$$
Note that $n(K)=n(K')+1$.  Since the vertex $u$ contributes to $q(K')$ but not to $q(K)$ and the vertex $v$
does not contribute to $q(K)$, we also have $q(K')=q(K)+1$.

Furthermore, we claim that $s(K')\ge s(K)$. Indeed, let $p_1$ and $p_2$ be the ends of $P$.  For $i\in\{1,2\}$, note
that Lemma~\ref{lemma:neartri} implies that $p_i$ has degree at least three in $G$, and thus degree at least two in $G-v$.
Let $u_i$ be the neighbor of $p_i$ in $\partial G-E(P)$ and let $w_i$ be the neighbor of $p_i$ in $\partial (G-v)-E(P)$.
If $p_iw_i$ does not contribute to $s(K')$, then we have $f(w_i)\le 3$, and thus $w_i$ is incident with the outer face of $G$,
and Lemma~\ref{lemma:1chords} implies that $w_i=u_i$.  We conclude that the edge $p_iw_i$ contributes at least as much to
$s(K')$ as the edge $p_iu_i$ contributes to $s(K)$.  Since only the two edges of $\partial G-E(P)$ incident with the
ends of $P$ can contribute to $s(K)$, we get the desired inequality $s(K')\ge s(K)$.

Therefore, we have
\begin{align*}
\lambda(K)&=n(K)+2q(K)+2s(K)+2\cdot(3-v(P))\\
&\le n(K')+1+2q(K')-2+2s'(K)+2\cdot(3-v(P))\\
&=\lambda(K')-1.
\end{align*}
This contradicts the inequality $\lambda(K)>\lambda(K')$.
\end{proof}

Though we are not able to eliminate 2-chords entirely (as we did for 1-chords), the following lemma characterizes them. 

\begin{lemma}\label{lemma:2chords}
    Let $K=(G,P,f)$ be a minimal counterexample.
    If $Q=uvw$ is a 2-chord of $G$, then one of the following conditions holds:
    \begin{enumerate}[(i)]
        \item\label{2ch:middle} $v(P) = 3$ and one of $u$ and $w$ is the middle vertex of $P$, or
	\item\label{2ch:triangle} $uvw$ is a triangle bounding an internal face of $G$ and $uw\in E(\partial G)\setminus E(P)$, or
        \item\label{2ch:tripod} the 2-chord $Q$ separates $G$ into subgraphs $G_1$ and $G_2$ such that $P\subseteq G_1$ and
	      the canvas $K_2=(G_2,Q,f)$ is a tripod.
    \end{enumerate}
\end{lemma}
\begin{proof}
    Let $G_1$ and $G_2$ be the subgraphs into which the $2$-chord $Q$ separates $G$, with $v(P\cap G_1)\ge v(P\cap G_2)$.
    If $P\not\subseteq G_1$, then \ref{2ch:middle} holds.  Hence, we can assume that $P\subseteq G_1$.
    Let $K_1=(G_1,P,f)$ and $K_2=(G_2,Q,f)$.  
    
    Note that the canvas $K_2$ is not \ref{S23a}-singular, as otherwise $K$ would have a vertex contradicting
    Lemma~\ref{lemma:minfvalues-notA}.  Since $v(Q)=3$, none of the conclusions of Observation~\ref{obs:comb} is possible,
    and thus either $Q$ is not an induced path in $G_2$, or $K_2$ is not tripod-free.
    Lemmas~\ref{lemma:noseptriang} and \ref{lemma:1chords} imply that \ref{2ch:triangle} holds in the former case
    and that \ref{2ch:tripod} holds in the latter case.
\end{proof}

Next we show that vertices of a minimal counterexample $(G,P,f)$ cannot have too many neighbors in $P$.

\begin{lemma}\label{lemma:nottoomanypnbrs}
    Let $K=(G,P,f)$ be a minimal counterexample with $v(P)\ge 2$.  If a vertex $v\in V(G)\setminus V(P)$ is
    adjacent to both ends of $P$, then $v(P)=3$, $v(\partial G)=4$, and $v$ is the unique vertex in $V(\partial G)\setminus V(P)$.
\end{lemma}
\begin{proof}
    Let $p$ and $p'$ be the endpoints of $P$ and let $Q=pvp'$.  Suppose first that $v\in V(\partial G)$.
    Then Lemma~\ref{lemma:1chords} implies that the cycle $\partial G$ is the concatenation of the paths $P$ and $Q$.
    If $v(P)=2$, then Lemma~\ref{lemma:noseptriang} implies that $\deg v=2$, contradicting Lemma~\ref{lemma:neartri}.
    Therefore, we have $v(P)=3$, and the outcome of the Lemma holds.
   
    Suppose now for a contradiction that $v\not\in V(\partial G)$, and thus $Q$ is a $2$-chord separating $G$ into
    subgraphs $G_1$ and $G_2$, where $P\subseteq G_1$.  Note that since $P$ is an induced path in $G$, the path $Q$ is
    induced in $G_2$.  Lemma~\ref{lemma:2chords} implies that $K_2=(G_2,Q,f)$ is a tripod.  That is, $V(G_2)\setminus V(Q)$
    consists of a single vertex $u$ adjacent to $p$, $p'$, and $v$ such that $f(u)=2$.

    Let $K'=(G-u,P,f')=\del(K,u)$; note that this operation is legal and has availability one.
    Since $v$ is the only vertex in $\partial (G-u)\setminus P$ and $f'(v)\ge 3$, the canvas $K'$ is valid and tripod-free.
    Moreover, an inspection of singular graphs shows that if $K'$ is singular, then $K'$ is \ref{S3}-singular.
    However, then $K$ would be \ref{S2i}-singular, which is a contradiction.
    Therefore, $K'$ is not singular.
    
    By the minimality of $K$, the canvas $K'$ is not a counterexample, and
    thus it is totally $c^{\lambda(K')}$-availably weakly degenerate.  It follows that $K$ is also
    totally $c^{\lambda(K')}$-availably weakly degenerate.  However, we have $n(K')=n(K)$,
    $q(K')\ge 0=q(K)$, and $s(K')\ge 0=s(K)$, and thus $\lambda(K')\ge \lambda(K)$.
    Therefore, $K$ is totally $c^{\lambda(K)}$-availably weakly degenerate, contradicting
    the assumption that $K$ is a counterexample.
\end{proof}

Finally, the following observation allows us to assume that all internal vertices have the $f$-value exactly four.

\begin{lemma}\label{lemma:internalfour}
Let $K=(G,P,f)$ be a minimal counterexample, and let $f'(v)=f(v)$ for every $v\in V(\partial G)\setminus V(P)$
and $f'(v)=4$ for every $v\in V(G)\setminus V(\partial G)$.  Then $K'=(G,P,f')$ is a minimal counterexample as well.
\end{lemma}
\begin{proof}
Since $K$ is valid and tripod-free, the canvas $K'$ is valid and tripod-free as well.
Moreover, since $K$ is non-singular, the canvas $K'$ is non-singular: This is clear if $V(G)\setminus V(\partial G)=\emptyset$,
as in that case $f'=f$.  Otherwise, an inspection of singular canvases excludes all cases except for $K'$ being \ref{S2i}-singular,
and if $K'$ were \ref{S2i}-singular, then $K$ would contradict Lemma~\ref{lemma:nottoomanypnbrs}.
Since $K$ is not totally $c^{\lambda(K')}$-availably weakly degenerate, Observation~\ref{obs:monot} shows that
$K'$ is not totally $c^{\lambda(K')}$-availably weakly degenerate, either.  Moreover, we clearly have $\lambda(K')=\lambda(K)$,
and thus $K'$ is a (minimal) counterexample.
\end{proof}

\section{Proof of Theorem \ref{thm:maintech}}\label{sec:proof}
Having established the basic structural properties of minimal counterexamples,
we are now equipped to prove Theorem \ref{thm:maintech}. The proof proceeds by
cases, considering all possible $f$-values for the first few vertices following $P$ in $\partial G$.
For a canvas $K=(G,P,f)$ such that $\partial G$ is a cycle, an \emph{initial segment}
of $K$ is a trail $pv_1\ldots v_k\subseteq \partial G$, where $p\in V(P)$ and $v_1\not\in V(P)$.

\begin{lemma}\label{lemma:notjustpv1}
    If $K=(G,P,f)$ is a minimal counterexample, then $v(\partial G) > v(P)+1$.
\end{lemma}
\begin{proof}
    Since $P$ is an induced path in $G$, we have $v(\partial G)\ge v(P)+1$.
    Suppose for a contradiction that $v(\partial G)=v(P)+1$, and thus $K$ has an initial segment $pvp'$
    with $p'\in V(P)$.  Note that this is only possible when $v(P)\ge 2$.  Let $U$ be the set of neighbors of $v$ different from $p$ and $p'$;
    Lemma~\ref{lemma:1chords} implies that $U\cap V(P)=\emptyset$, and Lemma~\ref{lemma:neartri} shows
    that $\partial (G-v)$ is a cycle with vertex set $V(P)\cup U$.  In particular we have $f_K(v)=f(v)-2\ge 0$.
    Moreover, $|U|\le 1$ is excluded by Lemmas~\ref{lemma:neartri} and \ref{lemma:nottoomanypnbrs}, and thus $|U|\ge 2$. 

    Let $K'=(G-v,P,f')=\del(K,v)$; this operation is legal and has availability $f_K(v)+1=f(v)-1\ge 1$.
    Since $f'(u)\ge 3$ for each $u\in U$, the canvas $K'$ is valid and tripod-free; and since $|U|\ge 2$, an inspection
    of singular canvases shows that it is non-singular.
    Hence, Observation~\ref{obs:mintolam} gives
    $$\lambda(K)>\lambda(K')+\log_c (f(v)-1).$$
    Note that Lemma~\ref{lemma:nottoomanypnbrs} implies that each vertex of $U$ has at most two neighbors in $P$,
    and thus $s(K')=2$.  Moreover, $n(K)\le n(K')+1$.  Consequently, we have
    \begin{align*}
    \lambda(K)&=n(K)+2q(K)+2s(K)+2\cdot(3-v(P))\\
    &\le n(K')+1+2q(K)+2s(K)+2\cdot(3-v(P))+2\cdot (s(K')-2)+2q(K')\\
    &=\lambda(K')+2q(K)+2s(K)-3,
    \end{align*}
    and thus
    $$2q(K)+2s(K)>3+\log_c (f(v)-1).$$
    However, if $f(v)=2$, then $q(K)=s(K)=0$; and otherwise $3+\log_c (f(v)-1)\ge \pgfmathparse{int(\CCC+3)}\pgfmathresult$
    and $q(K)+s(K)\le 3$.  In either case, we obtain a contradiction.
\end{proof}

\begin{lemma}\label{lemma:notjustpv1v2}
    If $K=(G,P,f)$ is a minimal counterexample, then $v(\partial G) > v(P)+2$.
\end{lemma}
\begin{proof}
    Suppose for a contradiction that $v(\partial G)=v(P)+1$, and thus $K$ has an initial segment $pv_1v_2p'$
    with $p'\in V(P)$.  Since $f^{-1}(2)$ is an independent set, we can by symmetry assume that $f(v_2)\ge 3$.
    
    Note that $p\neq p'$, as otherwise Lemma~\ref{lemma:noseptriang} implies that $\deg v_1=2$,
    contradicting Lemma~\ref{lemma:neartri}; and thus $v(P)\ge 2$.  Let $U$ be the set of neighbors of $v_1$
    different from $p$ and $v_2$; Lemma~\ref{lemma:1chords} implies that $U\cap V(\partial G)=\emptyset$, and
    by Lemma~\ref{lemma:neartri}, we have $V(\partial(G-v_1))=V(P)\cup U\cup\{v_2\}$
    and $|U|=\deg v_1-2\ge 1$.
    Lemma~\ref{lemma:1chords} also implies that $v_1$ and $v_2$ each have only one neighbor in $P$, and
    thus $f_K(v_1)=f(v_1)-1$ and $f_K(v_2)=f(v_2)-1$.

    Let $K'=(G',P,f')=\del(K,v_1)$; this operation is legal and has availability $f_K(v_1)+1=f(v_1)\ge 2$.
    Note that $f'(u)\ge 3$ for each $u\in U$ and $f'(v_2)=f(v_2)-1\ge 2$, and thus
    $K'$ is valid.  Moreover, Lemma~\ref{lemma:nottoomanypnbrs} implies that $K'$ is tripod-free.
    Since $|U|\ge 1$, an inspection of singular canvases shows that $K'$ is non-singular
    (let us remark that $K'$ is not \ref{S23a}-singular, since $v_2$ has only one neighbor in $P$,
    and it is not \ref{S23b}-singular, as otherwise the vertex in $U$ would violate Lemma~\ref{lemma:nottoomanypnbrs}).
    Hence, Observation~\ref{obs:mintolam} gives
    $$\lambda(K)>\lambda(K')+\log_c f(v_1).$$
    Note that $n(K')\ge n(K)-2$ (we remove $v_1$, and if $f'(v_2)=2$, then $v_2$ no longer contributes to the $n$-value).
    Moreover, note that $s(K')\ge 1$, since the vertex $u\in U$ that is adjacent to $p$ in $\partial(G')$ satisfies $f'(u)\ge 3$
    and has at most two neighbors in $P$ by Lemma~\ref{lemma:nottoomanypnbrs}; and that $s(K)\le 2$, since $\partial G$ is a cycle by Lemma~\ref{lemma:2conn},
    which implies that $s(K')\ge s(K)-1$.  Therefore,
    \begin{align*}
    \lambda(K)&=n(K)+2q(K)+2s(K)+2\cdot(3-v(P))\\
    &\le n'(K)+2+2q(K)+2\cdot(s(K')+1)+2\cdot(3-v(P))+2q(K')\\
    &=\lambda(K')+2q(K)+4.
    \end{align*}
    Consequently,
    $$2q(K)>(\log_c f(v_1))-4\ge \pgfmathparse{int(\CCC-4)}\pgfmathresult.$$
    Since $q(K)\le |V(\partial G)\setminus V(P)|=2$, this implies that $q(K)=2$, and thus $f(v_1)\ge 4$ and $f(v_2)\ge 4$.
    However, then we still get a contradiction, since in this case $(\log_c f(v_1))-4\ge\pgfmathparse{int(2*\CCC-4)}\pgfmathresult$.
\end{proof}

\begin{lemma}\label{lemma:v1leq3}
    If $pv$ is an initial segment of a minimal counterexample $K=(G,P,f)$, then $f(v)\le 3$.
\end{lemma}
\begin{proof}
    Suppose for a contradiction that $f(v)\ge 4$.  Let $f'(v)=f(v)-1$ and $f'(z)=f(z)$ for every $z\in V(G)\setminus (V(P)\cup\{v\})$,
    let $G'=G-pv$, and consider the canvas $K'=(G',P,f')$.  Note that $v(G')=v(G)$ and that the outer face of $G'$ is longer than
    the outer face of $G$, and thus $K'$ is smaller than $K$.  It follows that $K'$ is not a counterexample.

    The canvas $K'$ is clearly valid, and Lemma~\ref{lemma:nottoomanypnbrs} implies that it is tripod-free.
    By Lemma~\ref{lemma:neartri}, we have $\deg_G p\ge 3$, and thus $p$ has degree at least two in $G'$.  Lemma~\ref{lemma:1chords} implies that the neighbor $u$
    of $p$ in $\partial G'-E(P)$ is not incident with the outer face of $G$, and thus $f'(u)\ge 4$.
    It follows that $K'$ is not singular, and thus it is totally $c^{\lambda(K')}$-availably weakly degenerate.
    Observe that $f_K=f'_{K'}$, and thus $K$ is totally $c^{\lambda(K')}$-availably weakly degenerate as well.
    Moreover, observe that $\lambda(K')\ge \lambda(K)$---we have $n(K')=n(K)$, and moreover $q(K')\ge q(K)$ and $s(K')=s(K)$
    because of the vertex $u$.  Therefore, $K$ is totally $c^{\lambda(K)}$-availably weakly degenerate,
    which is a contradiction.
\end{proof}

\begin{lemma}\label{lemma:v2leq3}
    If $pv_1v_2$ is an initial segment of a minimal counterexample $K=(G,P,f)$, then $f(v_2)\le 3$.
\end{lemma}
\begin{proof}
    Suppose for a contradiction that $f(v_2)\ge 4$, and let $K'=(G-v_1,P,f')=\del(K,v_1)$.
    Lemma~\ref{lemma:1chords} implies that $p$ is the only neighbor of $v_1$ in $P$, and thus this
    operation has availability $f_K(v_1)+1=f(v_1)\ge 2$.  Let $U$ be the set of neighbors of $v_1$ in $G$
    different from $p$ and $v_2$; by Lemmas~\ref{lemma:1chords} and \ref{lemma:neartri},
    the set $U$ is non-empty and disjoint from $V(\partial G)$.  Note that $f'(v_2)\ge 3$
    and $f'(u)\ge 3$ for every $u\in U$.  This (together with the fact that $K$ is valid)
    implies that $K'$ is valid and non-singular, and Lemma~\ref{lemma:nottoomanypnbrs} implies that $K'$ is tripod-free.
    Hence, Observation~\ref{obs:mintolam} gives
    $$\lambda(K)>\lambda(K')+\log_c f(v_1)\ge\lambda(K')+\CCC.$$
    Note that $n(K')\ge n(K)-1$. Since $f(v_1)\le 3$ by Lemma~\ref{lemma:v1leq3}, all vertices of $K$ other than $v_2$ contributing to $q(K)$
    also contribute to $q'(K)$, and thus $q(K')\ge q(K)-1$.
    The edge $pu\in \partial(G-v_1)$ with $u\in U$ contributes $1$ to $s(K')$, since $u$ has at most two
    neighbors in $P$ by Lemma~\ref{lemma:nottoomanypnbrs}; and the edge of $\partial G-E(P)$ incident with the
    other end of $P$ contributes the same amount to $s(K)$ and to $s(K')$, since $v(\partial G)>v(P)+2$ by
    Lemma~\ref{lemma:notjustpv1v2}.  Therefore, we have $s(K')\ge s(K)$, and consequently
    \begin{align*}
    \lambda(K)&=n(K)+2q(K)+2s(K)+2\cdot(3-v(P))\\
    &\le n(K')+1+2\cdot(q(K')+1)+2s(K')+2\cdot(3-v(P))\\
    &=\lambda(K')+3.
    \end{align*}
    This is a contradiction.
\end{proof}

Recall that by Lemma \ref{lemma:notjustpv1v2}, each initial segment $pv_1v_2v_3$ of a minimal counterexample $(G,P,f)$
satisfies $v_3\not\in V(P)$.  Let us now restrict the $f$-value of $v_3$.
\begin{lemma}\label{lemma:v2notA,v3A}
    If $pv_1v_2v_3$ is an initial segment of a minimal counterexample $K=(G,P,f)$ and $f(v_2)\ge 3$, then $f(v_3)=2$.
\end{lemma}
\begin{proof}
    Suppose for a contradiction that $f(v_2)\ge 3$ and $f(v_3)\ge 3$, and let $K'=(G-v_1,P,f')=\del(K,v_1)$.
    Lemma~\ref{lemma:1chords} implies that $p$ is the only neighbor of $v_1$ in $P$, and thus this
    operation is legal and has availability $f_K(v_1)+1=f(v_1)\ge 2$.  Let $U$ be the set of neighbors of $v_1$ in $G$
    different from $p$ and $v_2$; by Lemmas~\ref{lemma:1chords} and \ref{lemma:neartri},
    the set $U$ is non-empty and disjoint from $V(\partial G)$.  Note that $f'(v_2)\ge 2$
    and $f'(u)\ge 3$ for every $u\in U$.  Moreover, the vertex $v_2$ has no neighbor $z\in V(G-v_1)\setminus P$
    with $f'(z)=2$, since otherwise the edge $v_2z$ would be a $1$-chord of $G$, contradicting Lemma~\ref{lemma:1chords}
    (note that $z\neq v_3$, since $f'(v_3)=f(v_3)\ge 3$).  Since $K$ is valid, this implies that $K'$ is valid as well.
    Moreover, Lemma~\ref{lemma:nottoomanypnbrs} implies that $K'$ is tripod-free.

    Note that $V(\partial (G-v_1))\setminus V(P)$ contains $v_3$ as well as a vertex $u\in U$ adjacent to $p$ in $\partial(G-v_1)$, and $f'(u)\ge 3$ and $f'(v_3)\ge 3$.  Therefore, the canvas $K'$ is not singular.
    Hence, Observation~\ref{obs:mintolam} gives
    $$\lambda(K)>\lambda(K')+\log_c f(v_1)\ge\lambda(K')+\CCC.$$
    Note that $n(K')\ge n(K)-2$ (the vertex $v_1$ was deleted and the vertex $v_2$ may not be contributing to $n(K')$).
    Since $f(v_1)\le 3$ and $f(v_2)\le 3$ by Lemmas~\ref{lemma:v1leq3} and \ref{lemma:v2leq3}, we have $q(K')=q(K)$.
    The edge $pu\in \partial(G-v_1)$ contributes $1$ to $s(K')$, since $u$ has at most two
    neighbors in $P$ by Lemma~\ref{lemma:nottoomanypnbrs}; and the edge of $\partial G-E(P)$ incident with the
    other end of $P$ contributes the same amount to $s(K)$ and to $s(K')$, since $v(\partial G)>v(P)+2$ by
    Lemma~\ref{lemma:notjustpv1v2}.  Therefore, we have $s(K')\ge s(K)$, and consequently
    \begin{align*}
    \lambda(K)&=n(K)+2q(K)+2s(K)+2\cdot(3-v(P))\\
    &\le n(K')+2+2q(K')+2s(K')+2\cdot(3-v(P))\\
    &=\lambda(K')+2.
    \end{align*}
    This is a contradiction.
\end{proof}

\begin{lemma}\label{lemma:v2inA,v3leq3}
    If $pv_1v_2v_3$ is an initial segment of a minimal counterexample $K=(G,P,f)$ and $f(v_2)=2$, then $f(v_3)=3$.
\end{lemma}
\begin{proof}
    Suppose for a contradiction that $f(v_2)=2$, but $f(v_3)\neq 3$.  Since the canvas $K$ is valid and $f(v_2)=2$, we have $f(v_3)>2$,
    and thus $f(v_3)\ge 4$.  Moreover, we also have $f(v_1)>2$, and Lemma~\ref{lemma:v1leq3} implies that $f(v_1)=3$.
    
    Let $K'=(G-v_2,P,f')=\del(K,v_2)$.  Lemma~\ref{lemma:1chords} implies that $v_2$ has no neighbors in $P$, and thus this
    operation is legal and has availability $f_K(v_2)+1=f(v_2)+1=3$.  Let $U$ be the set of neighbors of $v_2$ in $G$
    different from $v_1$ and $v_3$; by Lemmas~\ref{lemma:1chords} and \ref{lemma:neartri},
    the set $U$ is non-empty and disjoint from $V(\partial G)$.  Note that $f'(v_1)=2$
    and $f'(u)\ge 3$ for every $u\in U\cup \{v_3\}$.  Consequently, Lemma~\ref{lemma:1chords} implies that the canvas $K'$ is
    valid.  By Lemma~\ref{lemma:nottoomanypnbrs}, the canvas $K'$ is tripod-free.
    Moreover, since $|U\cup \{v_3\}|\ge 2$ and $U\cup \{v_3\}\subset V(\partial (G-v_2))\setminus V(P)$, the canvas $K'$ is non-singular.
    Hence, Observation~\ref{obs:mintolam} gives
    $$\lambda(K)>\lambda(K')+\log_c 3.$$
    Note that $n(K')=n(K)-1$ (the deleted vertex $v_2$ does not contribute to $n(K)$, and the vertex $v_1$ contributes to $n(K)$ but not to $n(K')$).
    Since $f(v_1)=3$, the only vertex that can contribute to $q(K)$ but not to $q(K')$ is $v_3$, and thus $q(K')\ge q(K)-1$.
    Similarly, we have $s(K')=s(K)-1$, since the only edge that contributes to $s(K)$ but not to $s(K')$ is $pv_1$. Therefore,
    \begin{align*}
    \lambda(K)&=n(K)+2q(K)+2s(K)+2\cdot(3-v(P))\\
    &\le n(K')+1+2\cdot(q(K')+1)+2\cdot(s(K')+1)+2\cdot(3-v(P))\\
    &=\lambda(K')+5.
    \end{align*}
    This is a contradiction.
\end{proof}

Lemmas~\ref{lemma:v1leq3} and \ref{lemma:v2leq3} imply that if $pv_1v_2v_3$ is an initial segment of a minimal counterexample $K=(G,P,f)$,
then $(f(v_1), f(v_2)) \in \{(3,3), (3,2), (2,3)\}$; and Lemmas~\ref{lemma:v2notA,v3A} and \ref{lemma:v2inA,v3leq3} furthermore give
$(f(v_1), f(v_2), f(v_3)) \in \{(3,3,2), (3,2,3), (2,3,2)\}$. The following lemmas rule out two of these three cases.

\begin{lemma}\label{lemma:not332}
   If $pv_1v_2v_3$ is an initial segment of a minimal counterexample $K=(G,P,f)$, then $(f(v_1), f(v_2), f(v_3)) \neq (3,3,2)$.
\end{lemma}
\begin{proof}
    Suppose for a contradiction that $(f(v_1), f(v_2), f(v_3))=(3,3,2)$.  Let $G'=G-pv_1$, $f'(v_1)=2$, and $f'(v)=f(v)$ for every $v\in V(G)\setminus (V(P)\cup \{v_1\})$,
    and consider the canvas $K'=(G',P,f')$.  
    By Lemma~\ref{lemma:neartri}, the graph $G$ has an internal face bounded by
    a triangle $pv_1u$, and since $\deg_G p\ge 3$, Lemma~\ref{lemma:1chords} implies that $u\not\in V(\partial G)$.  
    Since $v(G')=v(G)$ and the outer face of $G'$ is longer than the outer face
    of $G$, the canvas $K'$ is smaller than $K$, and thus it is not a
    counterexample.

    Note that Lemma~\ref{lemma:1chords} (together with the assumption that $f(v_2)>2$) implies that $K'$ is valid.
    Lemma~\ref{lemma:nottoomanypnbrs} implies that $K'$ is tripod-free.
    Since $f'(u)\ge 4$ and $u\in V(\partial G')$, the canvas $K'$ is non-singular.
    Since the canvas $K'$ is not a counterexample, it is totally $c^{\lambda(K')}$-availably weakly degenerate. Observe that $f_K=f'_{K'}$, and thus $K$ is totally $c^{\lambda(K')}$-availably weakly degenerate as well.
    Moreover, observe that $\lambda(K')>\lambda(K)$---we have $n(K')=n(K)-1$ because $v_1$ does not contribute to $n(K')$,
    $q(K')=q(K)+1$ since $u$ contributes to $q(K')$, and $s(K')\ge s(K)$ since the edge $pu$ contributing to $s(K')$ replaces the edge $pv_1$ contributing to $s(K)$.
    Therefore, $K$ is totally $c^{\lambda(K)}$-availably weakly degenerate, which is a contradiction.
\end{proof}

\begin{lemma}\label{lemma:not232}
    If $pv_1v_2v_3$ is an initial segment of a minimal counterexample $K=(G,P,f)$, then $(f(v_1), f(v_2), f(v_3))=(3,2,3)$.
\end{lemma}
\begin{proof}
  Suppose for a contradiction that $(f(v_1), f(v_2), f(v_3))\neq (3,2,3)$; as we have discussed before, this leaves the only possibility:
  $(f(v_1), f(v_2), f(v_3))=(2,3,2)$.  

  Let $K'=(G',P,f')=\del(\delsave(K, v_2, v_3), v_1)$; note that Lemma~\ref{lemma:1chords} implies that $f_K(v_2)=3$ and $f_K(v_1)=1$, and thus this sequence of operations
  is legal (though it can have availability one).  Let $u$ be the common neighbor of $v_1$ and $v_2$, which exists by Lemma~\ref{lemma:neartri}
  and is unique by Lemma~\ref{lemma:noseptriang}.  Let $Z$ be the set of vertices in $V(G)\setminus\{p,u,v_1,v_2,v_3\}$ with a neighbor in $\{v_1,v_2\}$.
  By Lemma~\ref{lemma:1chords}, the vertices of $Z\cup\{u\}$ are not contained in $V(\partial G)$,
  and by Lemma~\ref{lemma:internalfour}, we can assume that $f(z)=4$ for every $z\in Z\cup\{u\}$.
  Therefore, we have $f'(z)=3$ for $z\in Z$ and $f'(u)=2$.  Moreover, $f'(v_3)=f(v_3)=2$.

  Note that $u$ has no neighbor $w\in V(G')\setminus (V(P)\cup \{u\})$ such that $f'(w)=2$.  Indeed,
  $u$ is the only vertex of $K'$ such that $f(u)>2$ and $f'(u)=2$, and thus in that case we would have $f(w)=2$ as well;
  however, this would contradict Lemma~\ref{lemma:2chords} applied to the 2-chord $v_1uw$.  Since $K$ is valid, it follows that the canvas $K'$ is valid as well.
  Lemma~\ref{lemma:nottoomanypnbrs} implies that $K'$ is tripod-free.  
  
  Suppose next that the canvas $K'$ is singular.
  Lemma~\ref{lemma:minfvalues-notA} implies that $\deg_G v_2\ge 4$, and thus $|Z\cup\{u\}|\ge 2$, and $|V(\partial G')\setminus V(P)|\ge |Z\cup \{u,v_3\}|\ge 3$.
  An inspection of singular canvases shows that $v(P)=3$, $|Z\cup \{u,v_3\}|=3$, $V(G')=V(P)\cup Z\cup \{u,v_3\}$,
  and $v_3$ has a neighbor in $P$, i.e., $V(\partial G)\setminus V(P)=\{v_1,v_2,v_3\}$.
  It follows that $n(K)=3$ (only the two vertices of $Z\cup\{u\}$ and the vertex $v_2$ contribute) and $q(K)=s(K)=0$, and hence $\lambda(K)=3$.
  Moreover, since $v_3$ has a neighbor in $P$, we have $f_K(v_3)=1$, and thus the operation $\delsave(K, v_2, v_3)$ used in the process of obtaining $K'$
  has availability $2$.  Since $K'$ is weakly degenerate by Observation~\ref{obs:sing}, it follows that the canvas $K$ is totally $2$-availably weakly degenerate.
  Since $c^{\lambda(K)}=c^3<2$, this contradicts the assumption that $K$ is a counterexample.

  It follows that $K'$ is non-singular, and Observation~\ref{obs:mintolam} gives
  $$\lambda(K)>\lambda(K').$$
  Note that $n(K')=n(K)-2$ (the vertex $v_2$ is deleted, and the vertex $u$ does not contribute to $n(K')$) and $q(K')\ge q(K)$.
  Moreover, we have $s(K')\ge s(K)+s$, where $s=0$ if $u$ is adjacent to $P$ and $s=1$ otherwise: Indeed, Lemma~\ref{lemma:neartri} implies
  that the neighbor $z$ of $p$ in $\partial G'-E(P)$ belongs to $Z\cup\{u\}$; and if $z\in Z$, then the edge $pz$ contributes to $s(K')$,
  while the edge $pv_1$ does not contribute to $s(K)$.  Therefore,
  we have 
  \begin{align*}
  \lambda(K)&=n(K)+2q(K)+2s(K)+2\cdot(3-v(P))\\
  &\le n(K')+2+2q(K')+2\cdot(s(K')-s)+2\cdot(3-v(P))\\
  &=\lambda(K')+2-2s.
  \end{align*}
  It follows that $s=0$, i.e., that $u$ is the common neighbor of $p$, $v_1$, and $v_2$.
  Let us now distinguish two cases.

\paragraph{Case 1: $u$ is adjacent to at least two vertices of $P$.} By Lemma~\ref{lemma:nottoomanypnbrs}, it follows that $v(P)=3$, and,
letting $P=p_1p_2p$, the vertex $u$ is adjacent to $p_2$ and non-adjacent to $p_1$.  See Figure~\ref{fig:232case1}(a) for an illustration of this case,
and recall that by Lemma \ref{lemma:noseptriang}, the triangles $uv_1p$ and $upp_2$ are empty, implying that $\deg_G p=3$.
Let $G^\star=G'-p$, $P^\star=p_1p_2u$, and consider the canvas $K^\star = (G^\star, P^\star, f')$, as depicted in Figure~\ref{fig:232case1}(b).

 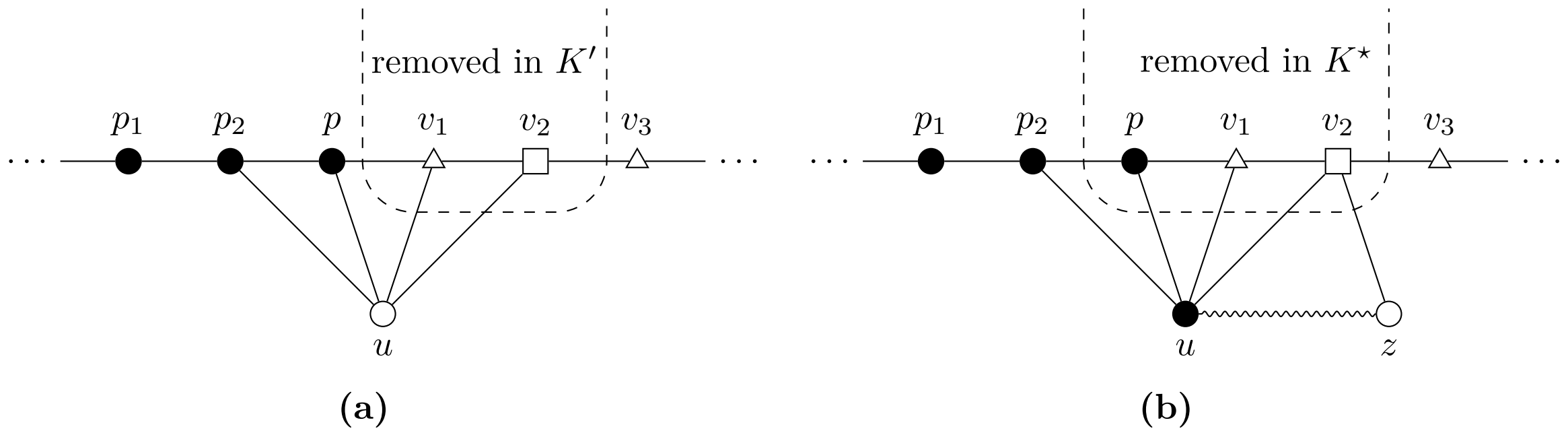
\begin{figure}[ht]\centering
    \subcaptionbox{}[0.44\textwidth]{\begin{tikzpicture}[scale=1]
    
    \node (ud) at (0,0){$\ldots$};
    \node[precoloured,label={above:$p_1$}] (p1) at (1,0){};
    
    \node[precoloured,label={above:$p_2$},fill=black] (p2) at (2,0){};

    \node[precoloured,label={above:$p$},fill=black] (p3) at (3,0){};
    
    \node[3vertex,label={above:$v_1$}] (v1) at (4,0){};
    
    \node[4vertex,label={above:$v_2$}] (v2) at (5,0){};

    \node[3vertex,label={above:$v_3$}] (v3) at (6,0){};
    
    \node (vd) at (7,0) {$\ldots$};

    \node[geq4vertex,label={below:$u$}] (u) at (3.5,-1.5){};

    \draw (ud) -- (p1) -- (p2) -- (p3) -- (v1) -- (v2) -- (v3) -- (vd);
    \draw (v1) -- (u) -- (v2);
    \draw (p3) -- (u) -- (p2);

    \draw[dashed] (3.3,1.5) -- (3.3,0) arc (180:270:0.5) -- (5.2,-0.5) arc (-90:0:0.5) -- (5.7,1.5);
        \node (R) at (4.5,1) {removed in $K'$};

\end{tikzpicture}

\begin{tikzpicture}[scale=1]
    
    \node (ud) at (0,0){$\ldots$};
    \node[precoloured,label={above:$p_1$}] (p1) at (1,0){};
    
    \node[precoloured,label={above:$p_2$},fill=black] (p2) at (2,0){};

    \node[precoloured,label={above:$p$},fill=black] (p3) at (3,0){};
    
    \node[3vertex,label={above:$v_1$}] (v1) at (4,0){};
    
    \node[4vertex,label={above:$v_2$}] (v2) at (5,0){};

    \node[3vertex,label={above:$v_3$}] (v3) at (6,0){};
    
    \node (vd) at (7,0) {$\ldots$};

    \node[geq4vertex,label={below:$u$},fill=black] (u) at (3.5,-1.5){};

    \node[geq4vertex,label={below:$z$}] (v) at (5.5,-1.5){};

    \draw (ud) -- (p1) -- (p2) -- (p3) -- (v1) -- (v2) -- (v3) -- (vd);
    \draw (v1) -- (u) -- (v2);
    \draw (p3) -- (u) -- (p2);
    \draw (v) -- (v2);
    \draw[wavy] (v) -- (u);

    \draw[dashed] (2.5,1.5) -- (2.5,0) arc (180:270:0.5) -- (5,-0.5) arc (-90:0:0.5) -- (5.5,1.5);
        \node (R) at (4.2,1) {removed in $K^\star$};

\end{tikzpicture}}%
    \hspace*{0.03\textwidth}
    \subcaptionbox{}[0.44\textwidth]{\begin{tikzpicture}[scale=1]
    
    \node (ud) at (0,0){$\ldots$};
    \node[precoloured,label={above:$p_1$}] (p1) at (1,0){};
    
    \node[precoloured,label={above:$p_2$},fill=black] (p2) at (2,0){};

    \node[precoloured,label={above:$p$},fill=black] (p3) at (3,0){};
    
    \node[3vertex,label={above:$v_1$}] (v1) at (4,0){};
    
    \node[4vertex,label={above:$v_2$}] (v2) at (5,0){};

    \node[3vertex,label={above:$v_3$}] (v3) at (6,0){};
    
    \node (vd) at (7,0) {$\ldots$};

    \node[geq4vertex,label={below:$u$}] (u) at (3.5,-1.5){};

    \draw (ud) -- (p1) -- (p2) -- (p3) -- (v1) -- (v2) -- (v3) -- (vd);
    \draw (v1) -- (u) -- (v2);
    \draw (p3) -- (u) -- (p2);

    \draw[dashed] (3.3,1.5) -- (3.3,0) arc (180:270:0.5) -- (5.2,-0.5) arc (-90:0:0.5) -- (5.7,1.5);
        \node (R) at (4.5,1) {removed in $K'$};

\end{tikzpicture}

\begin{tikzpicture}[scale=1]
    
    \node (ud) at (0,0){$\ldots$};
    \node[precoloured,label={above:$p_1$}] (p1) at (1,0){};
    
    \node[precoloured,label={above:$p_2$},fill=black] (p2) at (2,0){};

    \node[precoloured,label={above:$p$},fill=black] (p3) at (3,0){};
    
    \node[3vertex,label={above:$v_1$}] (v1) at (4,0){};
    
    \node[4vertex,label={above:$v_2$}] (v2) at (5,0){};

    \node[3vertex,label={above:$v_3$}] (v3) at (6,0){};
    
    \node (vd) at (7,0) {$\ldots$};

    \node[geq4vertex,label={below:$u$},fill=black] (u) at (3.5,-1.5){};

    \node[geq4vertex,label={below:$z$}] (v) at (5.5,-1.5){};

    \draw (ud) -- (p1) -- (p2) -- (p3) -- (v1) -- (v2) -- (v3) -- (vd);
    \draw (v1) -- (u) -- (v2);
    \draw (p3) -- (u) -- (p2);
    \draw (v) -- (v2);
    \draw[wavy] (v) -- (u);

    \draw[dashed] (2.5,1.5) -- (2.5,0) arc (180:270:0.5) -- (5,-0.5) arc (-90:0:0.5) -- (5.5,1.5);
        \node (R) at (4.2,1) {removed in $K^\star$};

\end{tikzpicture}}%
    \hspace*{0.01\textwidth}
    
    \captionsetup{width=0.9\textwidth}
    \caption{Case 1 of Lemma \ref{lemma:not232}. Vertices shown by triangles, squares, and empty circles have $f$-value $2$, $3$, and at least four,
    respectively. In (b), the wavy edge is counted by $s(K^\star)$. \label{fig:232case1}}
\end{figure}

Since $K'$ is valid and $p_1u\not\in E(G)$, the canvas $K^\star$ is valid as well.  Moreover, as we have argued above, the vertex $u$ has no neighbor $w\in V(G^\star)\setminus V(P^\star)$ with $f'(w)=2$.
This implies that $K^\star$ is tripod-free.  Recall also that $Z\cup\{v_3\}\subseteq V(\partial G^\star)\setminus V(P^\star)$ has size at least two. Since $G$ is a near-triangulation and $\deg v_2 \geq 4$ by Lemma~\ref{lemma:atleastfour}, $v_2$ and $u$ have a common neighbor $z$, where $z \not \in V(\partial G)$ by Lemma~\ref{lemma:1chords}.
An inspection of singular canvases shows that $K^\star$ is non-singular: Indeed, $|V(\partial G^\star)\setminus V(P^\star)|\ge 2$
and the fact that the end $u$ of $P^\star$ has no neighbor with $f'$-value two excludes all the cases except for \ref{S23a} and \ref{S23b}.
If $K^\star$ were \ref{S23a}-singular, then $v_3$ would be adjacent to $p_2$, contradicting Lemma~\ref{lemma:1chords}.
Finally, if $K^\star$ were \ref{S23b}-singular, then $K'$ would be \ref{S232b}-singular, which we have excluded before.

Consider the canvas $K''=(G'',P,f'')=\del(K',u)$.  This operation is legal with availability one, since $up_1\not\in E(G)$ implies $f'_{K'}(u)=0$.
Moreover, note that $G^\star-V(P^\star)=G''-V(P)$ and $f^\star_{K^\star}=f''_{K''}$.
Hence, Observation~\ref{obs:mintolam} (with $K''$ playing the role of $K'$) gives
$$\lambda(K)>\lambda(K^\star).$$
Note that $n(K^\star)=n(K')=n(K)-2$ and $q(K^\star)=q(K')=q(K)$.  Moreover, an edge $uz$ with $z\in Z$ belongs to $E(\partial G^\star)\setminus E(P)$ and contributes to $s(K^\star)$,
while the edge $pv_1$ does not contribute to $s(K)$, and thus we have $s(K^\star)\ge s(K)+1$.
However, this implies that 
\begin{align*}
\lambda(K)&=n(K)+2q(K)+2s(K)+2\cdot(3-v(P))\\
&\le n(K^\star)+2+2q(K^\star)+2\cdot(s(K^\star)-1)+2\cdot(3-v(P^\star))\\
&=\lambda(K^\star),
\end{align*}
which is a contradiction.

\paragraph{Case 2: The vertex $p$ is the only neighbor of $u$ in $P$.}  Thus, we have $f_K(u)=3$, while $f_K(v_1)=1$.
Note that Lemma~\ref{lemma:2chords} implies that $p$, $v_1$, and $v_2$ are the only neighbors of $u$ in $V(\partial G)$.
Let $z$ be the common neighbor of $u$ and $v_2$ different from $v_1$, which exists by Lemma~\ref{lemma:neartri},
and is unique by Lemma~\ref{lemma:noseptriang}.
Since $\deg_G u\ge 4$ by Lemma~\ref{lemma:neartri}, we have $z\neq v_3$, and thus $z\not\in V(\partial G)$ by Lemma~\ref{lemma:1chords}.
Moreover, note that Lemma~\ref{lemma:noseptriang} implies that $\deg_G v_1=3$.
Let us now distinguish two subcases.

\paragraph{Subcase 2a: $zv_3\not\in E(G)$.} Let
$$K^\dag=(G^\dag,P,f^\dag)=\Omit(\delsave(\delsave(K,v_2,v_3),u,v_1),v_1).$$
Observe that these operations are legal, and moreover, since $f_K(u)=f_K(v_1)+2$ and both $u$ and $v_1$ are adjacent to $v_2$,
the second $\delsave$ has availability two.
The canvas $K^\dag$ is tripod-free by Lemma~\ref{lemma:nottoomanypnbrs}.
Moreover, since $K$ is valid, $K^\dag$ is valid as well: Indeed, $z$ is the only vertex in $V(G^\dag)\setminus V(P)$ such that $f^\dag(z)=2$ but $f(z)>2$,
and since $zv_3\not\in E(G)$, Lemma~\ref{lemma:2chords} implies that $z$ has no neighbor $w\in V(G^\dag)\setminus (V(P)\cup\{z\})$ such that $f^\dag(w)=f(w)=2$
(otherwise, we would get a contradiction by considering the $2$-chord $v_2zw$, and noting that $wv_3\not\in E(G)$ since $f(v_3)=2$ and the canvas $K$ is valid).

Since $K^\dag$ has two distinct vertices $z$ and $v_3$ such that $f^\dag(z)=f^\dag(v_3)=2$,
it is easy to see that the canvas $K^\dag$ is not singular:  Indeed, otherwise it would have to be \ref{S232a}-, \ref{S232b}-, or \ref{S232c}-singular.
It cannot be \ref{S232a}-singular, since $v_3$ is not adjacent to the middle vertex of $P$ by Lemma~\ref{lemma:1chords},
and it cannot be \ref{S232c}-singular by Lemma~\ref{lemma:nottoomanypnbrs}.  Finally, if it were \ref{S232b}-singular, then either $v_3$ would be adjacent to the middle vertex of $P$,
contradicting Lemma~\ref{lemma:1chords}, or $z$ would be adjacent to $p$ and the $2$-chord $pzv_2$ would contradict Lemma~\ref{lemma:2chords}.
Hence, Observation~\ref{obs:mintolam} gives
$$\lambda(K)>\lambda(K^\dag)+\CCC.$$
However, note that $n(K^\dag)=n(K)-3$ (we have deleted $v_2$ and $u$, and $z$ does not contribute to $K^\dag$),
$q(K^\dag)=q(K)$, and $s(K^\dag)\ge s(K)$ (since the edge $pv_1$ does not contribute to $s(K))$,
and thus $\lambda(K^\dag)\ge \lambda(K)-3$.  This is a contradiction.

\paragraph{Subcase 2b: $zv_3\in E(G)$.}  By Lemma~\ref{lemma:noseptriang}, it follows that $v_1$, $u$, $z$, and $v_3$
are the only neighbors of $v_2$ in $G$, and thus $\deg_G v_2=4$.  Let $K^\ddag=(G^\ddag,P,f^\ddag)=\Omit(\Omit(\delsave(K,u,v_1),v_1),v_2)$.
These operations are legal, and the $\delsave$ has availability two.  As we have noted before, the vertex $u$ has no neighbors in $V(\partial G)\setminus\{p,v_1,v_2\}$,
and thus there is no vertex $w\in V(G^\ddag)\setminus V(P)$ with $f^\ddag(w)=2$ and $f(w)>2$.  Since the canvas $K$ is valid,
the canvas $K^\ddag$ is valid as well.  Moreover, Lemma~\ref{lemma:nottoomanypnbrs} implies that the canvas $K^\ddag$ is tripod-free.

Note that $z,v_3\in V(\partial G^\ddag)\setminus V(P)$, $f^\ddag(z)=3$, and $f^\ddag(v_3)=2$.
Moreover, Lemma~\ref{lemma:2chords} implies that $zp\not\in E(G)$ and that $z$ has no neighbor
$w\in V(\partial G^\ddag)\setminus (V(P)\cup \{z,v_3\})$ such that $f^\ddag(w)=f(w)=2$.
An inspection of the singular canvases shows that $K^\ddag$ is non-singular.
Hence, Observation~\ref{obs:mintolam} gives
$$\lambda(K)>\lambda(K^\ddag)+\CCC.$$
However, we have $n(K^\ddag)=n(K)-2$ (we have deleted $v_2$ and $u$),
$q(K^\ddag)=q(K)$, and $s(K^\ddag)\ge s(K)$ (since the edge $pv_1$ does not contribute to $s(K))$,
and thus $\lambda(K^\ddag)\ge \lambda(K)-2$.  This is a contradiction.
\end{proof}

The final case is unfortunately more complicated, and we need to consider the possible values of $f(v_4)$.
\begin{lemma}\label{lemma:atleastfour}
    If $K=(G,P,f)$ is a minimal counterexample, then $v(\partial G) > v(P)+3$, and moreover each initial segment
    $pv_1\dots v_4$ of $K$ satisfies $f(v_4)=2$.
\end{lemma}
\begin{proof}
    Suppose for a contradiction that this is not the case.  Hence, there exists an initial segment $pv_1v_2v_3$ of $K$
    (with $(f(v_1), f(v_2), f(v_3))=(3,2,3)$ by Lemma~\ref{lemma:not232}) such that $v_2$ is the only neighbor of $v_3$
    in $\partial G-V(P)$ with $f$-value $2$.  Lemma~\ref{lemma:1chords} then implies that $f(z)\ge 3$ holds for every neighbor
    $z\in V(G)\setminus(V(P)\cup \{v_2\})$ of $v_3$.

    Let $U$ be the set of neighbors of $v_2$ different from $v_1$ and $v_3$; Lemma~\ref{lemma:1chords} implies that
    $U\cap V(\partial G)=\emptyset$, and Lemma~\ref{lemma:minfvals-A} gives $|U|\ge 1$.
    Let $K'=(G',P,f')=\del(K,v_2)$.  This operation is legal and has availability 3, since $v_2$ has no neighbors in $P$ by Lemma~\ref{lemma:1chords}.
    Moreover, the only vertices in $V(G')\setminus V(P)$ whose $f'$-value is $2$ and $f$-value is greater than two are $v_1$ and $v_3$,
    and thus Lemma~\ref{lemma:1chords} together with the observation on $v_3$ from the previous paragraph imply that the canvas $K'$ is valid.
    It is also tripod-free by Lemma~\ref{lemma:nottoomanypnbrs}.
    Moreover, since $U\cup\{v_1,v_3\}\subseteq V(\partial G')\setminus V(P)$ and $U\neq\emptyset$, an inspection of
    the singular canvases shows that the canvas $K'$ is non-singular (note that $K'$ is neither \ref{S232a}- nor \ref{S232b}-singular,
    since Lemma~\ref{lemma:1chords} implies that neither $v_1$ nor $v_3$ is adjacent to the middle vertex of $P$;
    and it is not \ref{S232c}-singular by Lemma~\ref{lemma:nottoomanypnbrs}).
    Hence, Observation~\ref{obs:mintolam} gives
    $$\lambda(K)>\lambda(K')+\log_c 3.$$
    However, we have $n(K')=n(K)-2$ (the vertices $v_1$ and $v_3$ contribute to $n(K)$ and do not contribute to $n(K')$) and
    $q(K')=q(K)$.  Moreover, since $\partial G$ is a cycle, we have $s(K)\le 2$, and thus $s(K')\ge s(K)-2$.
    It follows that $\lambda(K')\ge \lambda(K)-6$, which is a contradiction.
\end{proof}

It remains to exclude the last possibility for the initial segment.

\begin{lemma}\label{lemma:finisher}
    No minimal counterexample exists.
\end{lemma}
\begin{proof}
Suppose for a contradiction that $K=(G,P,f)$ is a minimal counterexample.  By Lemma~\ref{lemma:internalfour}, we may assume that $f(z)=4$ for every $z\in V(G)\setminus V(P)$.
Lemma~\ref{lemma:atleastfour} implies that $v(\partial G) > v(P)+3$
and that an initial segment $pv_1\ldots v_4$ of $K$ satisfies $(f(v_1), f(v_2), f(v_3), f(v_4))=(3,2,3,2)$.
Let $K'=(G',P,f')=\del(\delsave(K,v_3,v_4),v_2)$.  Using Lemma~\ref{lemma:1chords}, we see that the operations in this sequence are legal,
and that the $\del$ operation has has availability 2.  Let $u$ be the common neighbor of $v_2$ and $v_3$ in $G$ (which exists by Lemma~\ref{lemma:neartri}, and is unique
by Lemma~\ref{lemma:noseptriang}), and for $i\in\{2,3\}$, let $U_i$ be the set of neighbors of $v_i$ in $V(G)\setminus \{v_{i-1},v_{i+1}, u\}$.
Note that Lemma~\ref{lemma:minfvalues-notA} implies that $U_3\neq \emptyset$.
Lemma~\ref{lemma:1chords} implies that $(\{u\}\cup U_2\cup U_3)\cap V(\partial G)=\emptyset$, and that each vertex $z\in V(G')\setminus V(P')$ satisfies $f'(z)=2<f(z)$
if and only if $z\in\{v_1,u\}$.  Note that $v_1$ has no neighbor $w\in V(G')\setminus V(P')$ such that $f(w)=2$ by Lemma~\ref{lemma:1chords},
and $u$ has no such neighbor by Lemma~\ref{lemma:2chords} applied to the hypothetical $2$-chord $v_2uw$ of $G$.

Suppose first that $uv_1\not\in E(G)$.  Then the preceding analysis shows that the canvas $K'$ is valid.
It is also tripod-free by Lemma~\ref{lemma:nottoomanypnbrs}, and non-singular since $|V(G')\setminus V(P')|\ge |\{v_1,u,v_4\}\cup U_3|\ge 4$.
Hence, Observation~\ref{obs:mintolam} gives
$$\lambda(K)>\lambda(K')+\CCC.$$
However, we have $n(K')=n(K)-3$ (the vertex $v_3$ is deleted and the vertices $v_1$ and $u$ do not contribute to $n(K')$),
$q(K')=q(K)$, and $s(K')\ge s(K)-1$ (the edge $pv_1$ contributes to $s(K)$ but not to $s(K')$, and the edge of $\partial G-E(P)$ incident
with the other end of $P$ is unaffected), and thus $\lambda(K')\ge \lambda(K)-5$.  This is a contradiction.

Therefore, we have $uv_1\in E(G)$.  The triangles $uv_1v_2$ and $uv_2v_3$ are empty by Lemma~\ref{lemma:noseptriang}, and thus
$\deg_G v_2=3$ and $U_2=\emptyset$.  Note that if $u$ has a neighbor $q\in V(P)$, then Lemma~\ref{lemma:2chords} applied to the $2$-chord $quv_2$
implies that $v(P)=3$ and $q$ is the middle vertex of $P$.
We now break into cases depending on whether $u$ has a neighbor in $V(P)$ or not.

\paragraph{Case 1: $v(P)=3$, $P=p_1p_2p$, and $u$ is adjacent to $p_2$.}
The 1-chord $p_2u$ of $G'$ separates $G'$ into two subgraphs $G_1$ and $G_2$, where $p_1\in V(G_1)$ and $p\in V(G_2)$.
Let $P_1=p_1p_2u$, $K_1=(G_1,P_1,f')$, $P_2=pp_2u$, and $K_2=(G_2,P_2,f')$.
Recall that we have shown above that $uv_1$ is the only edge of $G'-V(P)$ joining two vertices with $f'$-value $2$, and that $u$ is non-adjacent
to $p$ and $p_1$, and thus the canvases $K_1$ and $K_2$ are valid.  Lemma~\ref{lemma:1chords} implies that $p_2v_1\not\in E(G)$, and thus $K_2$
is tripod-free; and that all neighbors of $p_2$ in $V(\partial G_1)\setminus V(P_1)$ belong to $U_3$ (and have $f'$-value three), and thus $K_1$
is tripod-free.

For $i\in \{1,2\}$, let $\lambda_i=0$ if the canvas $K_i$ is singular and $\lambda_i=\lambda(K_i)$ if it is non-singular.
By Observation~\ref{obs:sing} and the minimality of $K$, we see the canvas $K_i$ is totally $c^{\lambda_i}$-availably weakly degenerate; let $S_i$ be an $\Omit$-free sequence of operations that witnesses this.
Consider the sequence of operations starting with the $\delsave$ and $\del$ (of availability two) used to create $K'$ from $K$,
continuing with $\del(K',u)$ (availability two), and finishing with the concatenation of $S_1$ and $S_2$.
All the operations in this sequence are legal and show that the canvas $K$ is totally $(4\cdot c^{\lambda_1+\lambda_2})$-availably weakly degenerate.
Since $K$ is a counterexample, it follows that
$$\lambda(K)>\lambda_1+\lambda_2+\pgfmathparse{int(2*\CCC)}\pgfmathresult.$$
Note that $n(K)=n(K')+3=n(K_1)+n(K_2)+3$ (since $u$ does not contribute to $n(K')$), $q(K)=q(K')=q(K_1)+q(K_2)$, and $s(K)\le s(K')+1\le s(K_1)+s(K_2)+1$
(if an edge $p_1v$ contributes to $s(K')$, it also contributes to $s(K_2)$, since $v\in V(\partial G)$ is not adjacent to $p_2$ by Lemma~\ref{lemma:1chords},
and thus $v$ has at most two neighbors in $V(P_2)$.  It follows that
\begin{align*}
\lambda(K)&=\lambda(K)-2\cdot(3-v(P))=n(K) + 2q(K) + 2s(K)\\
&\le(n(K_1)+n(K_2)+3)+ 2\cdot(q(K_1)+q(K_2)) + 2\cdot(s(K_1)+s(K_2)+1)\\
&=(\lambda(K_1)-2\cdot(3-v(P_1)))+(\lambda(K_2)-2\cdot(3-v(P_2)))+5\\
&=\lambda(K_1)+\lambda(K_2)+5.
\end{align*}
Therefore, we have
$$\lambda(K_1)+\lambda(K_2)>\lambda_1+\lambda_2+\pgfmathparse{int(2*\CCC-5)}\pgfmathresult,$$
and thus there exists $i\in\{1,2\}$ such that
$$\lambda(K_i)>\lambda_i+\pgfmathparse{(2*\CCC-5)/2}\pgfmathresult.$$
This is a contradiction, since $\lambda(K_i)=\lambda_i$ if $K_i$ is non-singular and $\lambda(K_i)\le 3$ by Observation~\ref{obs:sing} if $K_i$ is singular.

\paragraph{Case 2: $u$ has no neighbor in $V(P)$.}  Let $$K''=(G'',P,f'')=\Omit(v_2,\delsave(\delsave(K,v_3,v_4),u,v_2)).$$
Note that $f_K(u)=f(u)=4$ and $f_K(v_2)=f(v_2)=2$, and since both $u$ and $v_2$ are adjacent to $v_3$,
the second $\delsave$ is legal and has availability $2$; while the $\Omit$ operation is legal because of this $\delsave$
and because $\deg_G v_2=3=f_K(v_2)+1$.

Let $u'$ be the common neighbor of $u$ and $v_3$ different from $v_2$ (which exists by Lemma~\ref{lemma:neartri}, and is unique by Lemma~\ref{lemma:noseptriang}),
let $U'_3=U_3\setminus\{u'\}$, and let $U$ be the set of neighbors of $u$ in $V(G)\setminus (V(P)\cup\{v_1,v_2,v_3,u'\})$.
Note that $U\cap V(\partial G)=\emptyset$, as otherwise the $2$-chord $v_1uw$ for $w\in U\cap V(\partial G)$ would contradict Lemma~\ref{lemma:2chords}.
Hence, if a vertex $w\in V(G'')\setminus V(P)$ satisfies $f''(w)=2$ and $f(w)>2$, then $w\in\{v_1,u'\}$.
Moreover, Lemma~\ref{lemma:2chords} also implies that $u'$ has no neighbor $w\in V(G'')\setminus (V(P)\cup \{v_4\})$ satisfying $f''(w)=2$;
let us emphasize that this includes the case $w=v_1$, but not the case $w=v_4$.
Furthermore, Lemma~\ref{lemma:1chords} implies that $v_1$ has no neighbor $w\in V(G'')\setminus V(P)$ satisfying $f''(w)=2$.

If $u'v_4\not\in E(G)$, then the preceding analysis shows that the canvas $K''$ is valid.
It is also tripod-free by Lemma~\ref{lemma:nottoomanypnbrs}, and non-singular since there exist at least three vertices in $V(G'')\setminus V(P)$,
namely $v_1$, $u'$, and $v_4$, with $f''$-value two.
Hence, Observation~\ref{obs:mintolam} gives
$$\lambda(K)>\lambda(K'')+\CCC.$$
However, we have $n(K'')=n(K)-4$ (the vertices $u$ and $v_3$ are deleted and the vertices $v_1$ and $u'$ do not contribute to $n(K')$),
$q(K')=q(K)$, and $s(K')\ge s(K)-1$ (the edge $pv_1$ contributes to $s(K)$ but not to $s(K')$, and the edge of $\partial G-E(P)$ incident
with the other end of $P$ is unaffected), and thus $\lambda(K')\ge \lambda(K)-6$.  This is a contradiction.

     \begin{figure}[ht]\centering
{\begin{tikzpicture}[scale=1]
    
    \node (ud) at (2,0){$\ldots$};
    \node[precoloured,label={above:$p$},fill=black] (p3) at (3,0){};
    
    \node[4vertex,label={above:$v_1$}] (v1) at (4,0){};
    
    \node[3vertex,label={above:$v_2$}] (v2) at (5,0){};

    \node[4vertex,label={above:$v_3$}] (v3) at (6,0){};

    \node[3vertex,label={above:$v_4$}] (v4) at (7,0){};

    \node (vd) at (8,0) {$\ldots$};

    \node[geq4vertex,label={below:$u$}] (u) at (5,-1.5){};
    \node[geq4vertex,label={below:$u'$}] (u') at (6.2,-1.5){};

    \draw (ud) -- (p3) -- (v1) -- (v2) -- (v3) -- (v4) -- (vd);
    \draw (v1) -- (u) -- (v2);
    \draw (v3) -- (u);
    \draw (v3) -- (u') -- (v4);
    \draw (u) -- (u');

    \draw[dashed] (4.5,1.5) -- (4.5,-1.5) arc (180:270:0.6) -- (5.3,-2.1) arc (-90:0:0.6) -- (6.7,1.5);
        \node (R) at (5.5,1) {removed};

\end{tikzpicture}}%
    
    \captionsetup{width=0.9\textwidth}
    \caption{An illustration of the canvas $K$ in the final argument of Lemma~\ref{lemma:finisher}.\label{fig:lastclaimlastpart}}
\end{figure}
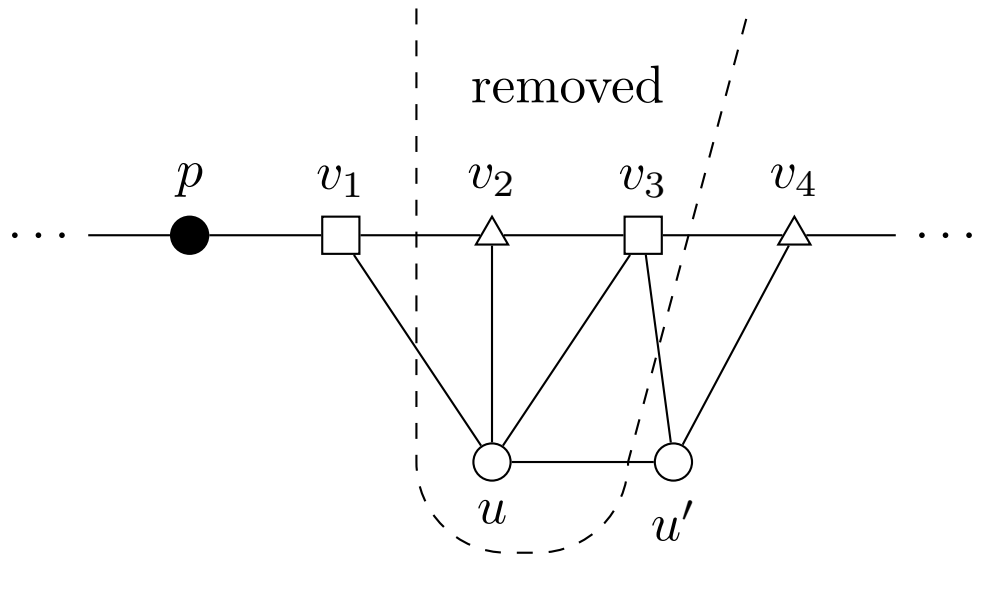

Therefore, we have $u'v_4\in E(G)$, and Lemma~\ref{lemma:noseptriang} implies that $\deg_G v_3=4$.
See Figure~\ref{fig:lastclaimlastpart} for an illustration of this case. 
Finally, let $K^\star=(G^\star, P,f^\star)=\Omit(\Omit(\delsave(K,u,v_2),v_2),v_3)$.
Since $f_K(u)=4$ and $f_K(v_2)=2$, the $\delsave$ operation has availability $2$,
and since $\deg_G v_2=3=f_K(v_2)+1$ and $\deg_G v_3=4=f_K(v_3)+1$,
the $\delsave$ makes the sequence of $\Omit$s legal.

The preceding analysis shows that there does not exist any vertex $w\in V(G^\star)\setminus (V(P)\cup\{v_1\})$ with $f^\star(w)=2$ and $f(w)>2$,
since $U\cap V(\partial G)=\emptyset$ and $f^\star(u')=3$.
Since $v_1$ has no neighbor $w\in V(G^\star)\setminus V(P)$ with $f(w)=2$ by Lemma~\ref{lemma:1chords},
it follows that the canvas $K^\star$ is valid.  It is also tripod-free by Lemma~\ref{lemma:nottoomanypnbrs},
and non-singular since $v_1,u',v_4\in V(G^\star)\setminus V(P)$ and $u'v_1\not\in E(G)$ (as we have observed before using Lemma~\ref{lemma:2chords}).
Hence, Observation~\ref{obs:mintolam} gives
$$\lambda(K)>\lambda(K^\star)+\CCC.$$
However, $n(K^\star)=n(K)-3$ (the vertices $u$ and $v_3$ are deleted, and the vertex $v_1$ does not contribute to $n(K^\star)$),
$q(K^\star)=q(K)$, and $s(K^\star)\ge s(K)-1$ (the edge $pv_1$ contributes to $s(K)$ but not to $s(K^\star)$, and the edge of $\partial G-E(P)$ incident
with the other end of $P$ is unaffected), and thus $\lambda(K^\star)\ge \lambda(K)-5$.  This is a contradiction.
\end{proof}

The proof of our main result is now trivial.

\begin{proof}[Proof of Theorem~\ref{thm:maintech}]
If Theorem~\ref{thm:maintech} were false, then there would exist a counterexample, and thus there would also exist a minimal one.
However, this is excluded by Lemma~\ref{lemma:finisher}.
\end{proof}

\section{AI Declaration}
We did not use AI at any point during the production of this manuscript.
\bibliographystyle{habbrv}
\bibliography{bibliog}

\noindent
\end{document}